\documentclass[12pt]{amsart}
\usepackage{amssymb,amsmath,amscd,verbatim}
\usepackage[all]{xy}
\font\tenrsf=rsfs10 at 11pt
\font\sevenrsf=rsfs7 at 8pt
\font\fiversf=rsfs5 at 6pt
\newfam\rsffam
\textfont\rsffam=\tenrsf
\scriptfont\rsffam=\sevenrsf
\scriptscriptfont\rsffam=\fiversf
\def\rond#1{{\tenrsf\fam\rsffam#1}}
\begin{document}
\newtheorem*{te*}{Theorem}
\newtheorem{cor}{Corollary}
\newtheorem{te}[cor]{Theorem}
\newtheorem{p}[cor]{Proposition}
\newtheorem*{lemma*}{Lemma}
\newtheorem{lemma}[cor]{Lemma}
\theoremstyle{definition}
\newtheorem*{defi*}{Definition}
\newtheorem{defi}[cor]{Definition}
\theoremstyle{remark}
\newtheorem{ob}[cor]{Remark}
\newtheorem{example}[cor]{Example}
\newcommand{\cC}{\mathcal{C}}
\newcommand{\cI}{\mathcal{I}}
\newcommand{\cO}{\mathcal{O}}
\newcommand{\cN}{\mathcal{N}}
\newcommand{\cS}{\mathcal{S}}
\newcommand{\cT}{\mathcal{T}}
\newcommand{\cV}{\mathcal{V}}
\newcommand{\dA}{d_A}
\newcommand{\Dr}{\rond{D}}
\newcommand{\Gr}{\rond{G}}
\newcommand{\Hr}{\rond{H}} 
\newcommand{\Kr}{\rond{K}}
\newcommand{\Or}{\rond{O}}
\newcommand{\Qr}{\rond{Q}}
\newcommand{\Bc}{\mathcal{B}}
\newcommand{\Dc}{\mathcal{D}}
\newcommand{\Kc}{\mathcal{K}}
\newcommand{\Oc}{\mathcal{O}}
\newcommand{\ccV}{{}^c\cV}
\newcommand{\ctX}{{}^cT\oX}
\newcommand{\ctsX}{{}^cT^*\oX}
\newcommand{\cun}{\cC^{\infty}}
\newcommand{\cunc}{\cC^{\infty}_c}
\newcommand{\cz}{\mathbb{C}}
\newcommand{\Diff}{\mathrm{Diff}}
\newcommand{\Diffc}{\mathrm{Diff}_c}
\newcommand{\Dom}{\mathrm{Dom}}
\newcommand{\dtt}{{d_{\tit}}} 
\newcommand{\End}{\mathrm{End}}
\newcommand{\ess}{\mathrm{ess}}
\newcommand{\Hom}{\mathrm{Hom}}
\newcommand{\im}{\mathrm{Im}} 
\newcommand{\Krke}{\Kr_{k,\varepsilon}}
\newcommand{\loc}{\mathrm{loc}}
\newcommand{\oX}{\overline{X}}
\newcommand{\px}{\partial_x}
\newcommand{\rz}{\mathbb{R}}
\newcommand{\supp}{\mathrm{supp}}
\newcommand{\C}{\mathbb{C}}
\newcommand{\R}{\mathbb{R}}
\newcommand{\N}{\mathbb{N}}
\newcommand{\Z}{\mathbb{Z}}
\newcommand{\siges}{\sigma_\ess} 
\newcommand{\Spec}{\mathrm{Spec}}
\newcommand{\ta}{\tilde{a}}
\newcommand{\tal}{\tilde{\alpha}}
\newcommand{\tih}{\tilde{h}}
\newcommand{\tf}{\tilde{\varphi}}
\newcommand{\tif}{\tilde{\phi}}
\newcommand{\tit}{\tilde{\theta}}
\newcommand{\tD}{\tilde{\Delta}}
\newcommand{\tM}{\tilde{M}}
\newcommand{\tmu}{\tilde{\mu}}
\newcommand{\ucz}{\underline{\cz}}
\newcommand{\Vol}{\operatorname{Vol}}
\newcommand{\Xe}{X_\varepsilon}
\newcommand{\zz}{\mathbb{Z}}
\def\cchi{\raisebox{.45 ex}{$\chi$}}
\def\l{<} 
\def\g{>} 
 
\title[Magnetic Schr\"odinger operators]{The spectrum of
magnetic Schr\"odinger operators and $k$-form Laplacians on conformally 
cusp manifolds} 
 
\author{Sylvain Gol\'enia} 
\address{Institutul de Matematic\u{a} al Academiei Rom\^{a}ne\\ 
P.O. Box 1-764\\RO-014700 
Bucha\-rest, Romania} 
\email{sylvain.golenia@imar.ro} 
\author{Sergiu Moroianu} 
\thanks{The authors were partially supported by the
contract MERG 006375, funded through the European Commission} 
\address{Institutul de Matematic\u{a} al Academiei Rom\^{a}ne\\ 
P.O. Box 1-764\\RO-014700 
Bucha\-rest, Romania} 
\email{moroianu@alum.mit.edu} 
\date{\today} 
\begin{abstract} 
We consider open manifolds which are interiors of a compact manifold
with boundary, and Riemannian metrics asymptotic to a conformally cylindrical
metric near the boundary. We show that the essential spectrum of the Laplace
operator on functions vanishes under the presence of a
magnetic field which does not define an integral relative cohomology class.
It follows that the essential spectrum is not stable by perturbation even by
a compactly supported magnetic field. We also treat magnetic operators
perturbed
with electric fields. In the same context we describe
the essential spectrum of the $k$-form Laplacian. This is shown to vanish
precisely when the $k$ and $k-1$ de Rham cohomology groups
of the boundary vanish.
In all the cases when we have pure-point spectrum we give Weyl-type
asymptotics for the eigenvalue-counting function. In the other cases 
we describe the essential spectrum.
\end{abstract} 
\maketitle 
 
\section*{Introduction} 
There exist complete, noncompact manifolds on which the
scalar Laplacian has pure point spectrum (e.g., \cite{donelili}).
One of the goals of this paper is to exhibit the same phenomenon 
for the Laplacian on differential forms. We study Riemannian manifolds $X$ 
with ends diffeomorphic to a cylinder $[0,\infty)\times M$,
where $M$ is a closed, possibly disconnected Riemannian manifold. The metric
on $X$ is the asymptotically conformally cylindrical
metric already studied in \cite{wlom}.
The general form of such a metric is given in Eq.\ \eqref{pme}; for the
purpose of this introduction, the reader should bear in mind the toy example
\begin{align}\label{mc} 
g_p=y^{-2p}(dy^2+h),&&y\to\infty\end{align} 
where $h$ is a metric on $M$ and $p>0$.
We first give conditions for the Laplacian on
$\Lambda^k X$ to have pure point spectrum. 
 
\begin{te}\label{th:0} 
Let $p>0$, $g_p$ the metric given by \eqref{pme} and 
\[\Delta_{k,p}:\cunc(X,\Lambda^kX)\to L^2(X,\Lambda^kX)\] 
the Laplacian on $k$-forms on $X$, viewed as an unbounded operator 
in $L^2$. Assume that the Betti numbers $h^k(M)$ and $h^{k-1}(M)$ 
both vanish. Then $\Delta_{k,p}$ is essentially self-adjoint
and has pure point spectrum.
\end{te} 
 
This result will be restated in larger generality in Theorem \ref{th3}. 
Note that $X$ is complete
if and only if $p\leq 1$. In that case, it is well-known that 
the Laplacian on forms is essentially self-adjoint. 
Without the hypothesis on the cohomology of the end, we can describe
the essential spectrum of $\Delta_{k,p}$. 
 
\begin{te}\label{th1} 
Assume that not both $h^k(M)$ and $h^{k-1}(M)$ vanish.
Let $g_0$ be an \emph{exact} cusp metric (see Definition
\ref{defex}), $0<p\leq 1$ and $g_p=x^{2p} g_0$. Then 
the essential spectrum of $\Delta_{k,p}$ is that of a direct sum
of at most two ordinary differential equations on the
real half-line, and
\begin{itemize}
\item If $0<p<1$ then $\siges(\Delta_{k,p})=[0,\infty)$. 
\item If $p=1$ then $\siges(\Delta_{k,1})=[c,\infty)$ where 
\[c=\begin{cases} 
\left(\frac{n-2k-1}{2}\right)^2&\text{if $h^{k-1}(M)=0$;}\\ 
\left(\frac{n-2k+1}{2}\right)^2&\text{if $h^{k}(M)=0$;}\\ 
\min\left\{\left(\frac{n-2k-1}{2}\right)^2,\left(\frac{n-2k+1}{2}\right)^2
\right\}&\text{otherwise.}\end{cases}\] 
\end{itemize} 
\end{te} 
 
In particular, the Laplacian of $g_p$ on $0$-, $1$- and $n$-forms always has 
non-empty essential spectrum when $X$ is complete, i.e., $p\leq 1$.

For $p>1$ we have a partial result, namely we show that all self-adjoint
extensions of $\Delta_{k,p}$ have pure point spectrum but only for the metric
\eqref{mc} (see Proposition \ref{p:essDelta}).

In a recent paper Antoci \cite{francesca1} 
studies the essential spectrum of a manifold $X$ with 
ends as above, with metric on the end 
\begin{equation}\label{anto} 
e^{-2(a+1)t}dt^2+e^{-2bt}h. 
\end{equation} 
Antoci assumes that $X$ is complete so $a\leq -1$. But with the change of 
variables $y=e^{(b-a-1)t}$ the metric \eqref{anto} 
is a constant multiple of \eqref{mc} under the assumption that 
$b-a-1> 0$ and $b <0$. Thus our results have some non-empty intersection 
with \cite{francesca1}; we improve one of Antoci's results by showing
that $0$ cannot be isolated in the essential spectrum for a general $M$.
(for $M=S^{n-1}$ this is proved in \cite[Theorem 6.1]{francesca1}).
This issue is actually important, since 
the absence of $0$ from the essential spectrum is equivalent to the
Hodge decomposition
\[L^2(X,\Lambda^k X,g_p)=\ker(\Delta_{k,p})\oplus\im(\Delta_{k,p}).\]
We also get somewhat more
refined information on the nature of the essential spectrum. Namely for
the metric \eqref{mc} we show the absence of singular continuous spectrum
(a classical fact for $k=0$).

The second topic of this article is the so-called magnetic
Laplacian on a manifold $X$ as above. We develop this in
Section \ref{sectionmagn}.
A magnetic field $B$ is a smooth exact $2$-form on $X$.
There exists a real $1$-form $A$, called magnetic potential,
satisfying $dA=B$. Set $\dA:=d+iA\wedge: \cun(X)\to\cun(X,T^*X)$.
The magnetic Laplacian on $\cun(X)$ is given by
\begin{equation}\label{eq:intro1}
\Delta_A:=\dA^*\dA
\end{equation}
When the manifold is complete,
this operator is essentially self-adjoint \cite{shubin}.
Given two magnetic potentials $A$ and $A'$ such that $A-A'$ is exact,
the two magnetic Laplacians $\Delta_A$ and $\Delta_{A'}$
are unitarily equivalent (this property is called gauge invariance).
Hence when $H^1(X,\rz)=0$, the spectral properties of the magnetic Laplacian do
not depend on the choice of the magnetic potential $A$.

We have found a strong relationship between the emptiness of the essential
spectrum of $\Delta_A$ and the cohomology class $[B]$
inside the relative cohomology group $H^2(X,M)$. We show first
that if the magnetic field does not vanish on the boundary, 
then the operator $\Delta_A$ has compact resolvent. 
In comparison, in
$\rz^n$ with the Euclidean metric this happens when the magnetic field
blows up at infinity with a certain control on its derivatives (see \cite{hm}
and references therein). Furthermore, if $B$
does vanish on $M$ but it defines a non-integral cohomology class $[B]$
inside the relative cohomology group $H^2(X,M)$ (up to a factor of $2\pi$), 
then the same conclusion holds. 

In $\rz^n$, 
a compactly supported magnetic field does not affect the essential 
spectrum \cite{CFKS, MP} and furthermore the wave operators are asymptotically 
complete for a certain choice of gauge \cite{LT}.
In contrast, on a manifold $X$ as above,
in dimensions other than $3$ we find in particular that a generic magnetic
field $B$ with compact support will turn off the essential spectrum
of $\Delta$. 
It seems quite surprising on one hand that we get the (non)integrality
condition, and on the other hand that there are no orientable
examples in dimension $3$. We stress that it is not the size of the magnetic 
field which kills the essential spectrum, but a more subtle 
phenomenon which has yet to be interpreted physically. 

In the case of $H^1(X,\rz)\neq 0$, in quantum mechanics, it is known that 
the choice of a potential vector has a physical meaning. This is known 
as the Aharonov-Bohm effect \cite{AB}. In $\rz^n$ with a bounded
obstacle, this  phenomenon can be seen through a difference of phase
of the waves arising from two non-homotopic paths that circumvent the  
obstacle. Note that the essential spectrum is independent of this 
choice. In contrast, in our setting we show that the essential spectrum itself 
can be turned off by a suitable choice of magnetic potential.
 
In all these instances when the operators have empty
essential spectra, the eigenvalues turn out to satisfy an asymptotic
law that we prove using results from \cite{wlom}. We state here only one
such result and we refer to the body of the paper for the others.
 
\begin{te}\label{th3} 
In the setting of Theorem \ref{th:0}, assume that $h^k(M)=h^{k-1}(M)=0$. 
Let $N_{k,p}$ be the counting function of the
eigenvalues of $\Delta_{k,p}$. Then 
\begin{equation}\label{rg} 
N_{k,p}(\lambda) \approx \begin{cases} 
C_1\lambda^{n/2}& 
\text{for $1/n< p<\infty$,}\\
C_2\lambda^{n/2}\log \lambda &\text{for $p=1/n$,}\\ 
C_3\lambda^{1/2p}&\text{for $0<p<1/n$} 
\end{cases} 
\end{equation} 
where
\begin{equation}\label{c12}\begin{split} 
C_1=&{\dbinom{n}{k}}\frac{\Vol(X,g_p)\Vol(S^{n-1})}{n(2\pi)^n},\\ 
C_2=&{\dbinom{n}{k}} \frac{\Vol(M,{a_0}^{1/2}dh_0)\Vol(S^{n-1})}{2(2\pi)^n}.
\end{split}\end{equation} 
If moreover we assume that $g_0$ is an exact cusp metric, then 
\[C_3=\frac{\Gamma\left(\tfrac{1-p}{2p}\right) 
\left(\zeta\left(\Delta_k^{h_0}, \frac{1}{p}-1\right) 
+\zeta\left(\Delta_{k-1}^{h_0}, \frac{1}{p}-1\right)\right)} 
{2\sqrt{\pi}\Gamma\left(\tfrac{1}{2p}\right)},\]
where $\Delta_k^{h_0}$ is the Laplacian on $k$ forms on $M$ with
respect to the canonical metric $h_0$ defined in Section \ref{lcm}. 
\end{te}
 
We also find some surprising conditions for the Schr\"o\-din\-ger operator
on functions to have pure-point spectrum.
It is well-known that if the electric potential blows up at infinity
then there is no essential spectrum. But we find potentials
which are not bounded below for which the same happens.

The result continues to hold in the more general setting of magnetic
Schr\"o\-din\-ger operators acting on $k$-forms, for which 
we also compute the eigenvalue asymptotics. When the magnetic field is 
non-integral, the potential (which can be endomorphism-valued) may even 
go uniformly to $-\infty$ at a specified rate. Finally,
even if the potential is not self-adjoint, we give some conditions 
for the absence of the essential spectrum.

Our methods are based on the Melrose cusp calculus of pseudodifferential 
operators (see e.g., \cite{meni96c}). It would also be interesting to 
apply in this setting the groupoid techniques of \cite{lani}.

{\bf \small Acknowledgments: } 
{\small We thank Barbu Berceanu, Dan Burghelea, Vladimir Georgescu, 
Marius M\u antoiu and Radu Purice for helpful discussions.} 
 
\section{The Laplacian of a cusp metric}\label{lcm} 
 
This section follows closely \cite{wlom}, see also \cite{meni96c}.
Let $\oX$ be a smooth $n$-dimensional compact manifold
with closed boundary $M$, and $x:\oX\to[0,\infty)$ a
boundary-defining function. A \emph{cusp metric} on $\oX$ is a
complete Riemannian metric $g_0$ on $X:=\oX\setminus M$
which in local coordinates near the boundary takes the form 
\begin{equation}\label{cume} 
g_0=a_{00}(x,y)\frac{dx^2}{x^4}+\sum_{j=1}^{n-1} a_{0j}(x,y) 
\frac{dx}{x^2}dy_j+\sum_{i,j=1}^{n-1}a_{ij}(x,y)dy_idy_j 
\end{equation} 
such that the matrix $(a_{\alpha\beta})$ is smooth and non-degenerate 
down to $x=0$. For example, if $a_{00}=1$, $a_{0j}=0$ and
$a_{ij}$ is independent of $x$, we get a product metric near $M$. 
Thus a cusp metric is asymptotically cylindrical, with a 
precise meaning for the word ``asymptotic". 
 
Let $\cI\subset\cun(\oX)$ be the principal ideal generated by the function
$x$. Recall that a \emph{cusp vector field} is a smooth vector field $V$ on
$\oX$ such that $dx(V)\in\cI^2$. The space of cusp vector fields forms a
Lie subalgebra $\ccV$ of the Lie algebra of smooth vector fields on $\oX$.
Let $E,F\to\oX$ be smooth vector bundles. The space of cusp differential
operators $\Diffc(\oX,E,F)$ is the space of those differential operators
which in local trivializations can be written as composition of cusp
vector fields and smooth bundle morphisms. 
 
The \emph{normal operator} of $P\in\Diffc(\oX,E,F)$ is defined by 
\[\rz\ni\xi\mapsto\cN(P)(\xi):= 
\left(e^{i\xi/x}Pe^{-i\xi/x}\right)_{|x=0}\in\Diff(M,E_{|M},F_{|M}).\] 
From the definition, $\ker \cN=\cI\cdot\Diffc$, which we denote again by $\cI$. 
 
We want to study the Laplacian on $k$-forms associated to the metric 
\begin{equation}\label{pme}
g_p:=x^{2p}g_0
\end{equation}
where $g_0$ is a cusp metric given near $M$ by \eqref{cume}. Note that the
metric \eqref{mc} is a particular case of such metric (set $x=1/y$).
By a general principle, $\Delta_{k,p}$ should be a cusp differential 
operator of type $(2,2p)$, i.e., $\Delta_{k,p}\in
x^{-2p}\Diffc^2(\oX,\Lambda^k(\oX))$. We show below that this is so, and we
compute the normal operator of $x^{2p}\Delta_{k,p}$. 
 
In \cite[Lemma 6]{wlom} it is constructed a canonical metric on
$M$ induced from $g_0$. From now on we fix a product decomposition 
of $\oX$ near $M$ and we rewrite \eqref{cume} as
\begin{equation}\label{cume2}
g_0=a \left(\frac{dx}{x^2}+\alpha(x)\right)^2+h(x)
\end{equation}
where
\begin{align*} 
a:=a_{00}&\in\cun(\oX)\\ 
\alpha&\in\cun([0,\varepsilon)\times M,\Lambda^1(M))\\ 
h&\in \cun([0,\varepsilon)\times M, S^2TM). 
\end{align*} 
Then by \cite[Lemma 6]{wlom}, the function $a_0:=a(0)$, 
the metric $h_0:=h(0)$ and the class modulo exact forms 
of the $1$-form $\alpha_0:=\alpha(0)$, defined on $M$,
are independent of the chosen product decomposition and of 
the boundary-defining function $x$ inside the fixed cusp structure.

\begin{defi}\label{defex}
The metric $g_0$ is called \emph{exact} if $a_0=1$ and $\alpha_0$ is an
exact $1$-form.
\end{defi}
If $g_0$ is exact, by replacing $x$ with another boundary-defining
function inside the same cusp structure, we can as well assume that 
$\alpha=0$ (see \cite{wlom}).

Let $V_0$ be the vector field $x^2\px\in\ccV$ (this makes
sense since we fixed a product decomposition). Choose local coordinates 
$(y_i)_{1\leq i\leq n-1}$ on $M$ and set
\begin{align}\label{vi} 
V_i:=&\partial_{y_i}-\alpha(\partial_{y_i})V_0,&i=1,\ldots,n-1. 
\end{align}
Notice that $V_0$ is orthogonal to $V_i$, $i\geq 1$ with respect to $g_0$.
Let
\begin{align*} 
V^0:=&x^{-2}dx+\alpha,& V^i=&dy^i 
\end{align*}
be the dual basis. Then $(V^j)_{0\leq j\leq n-1}$ form a local basis
for $\ctsX$ near $M$. Thus we get an orthogonal decomposition of
smooth vector bundles 
\begin{equation}\label{decom} 
\Lambda^k(\ctX)\simeq\Lambda^k(TM)\oplus V^0\wedge\Lambda^{k-1}(TM) 
\end{equation} 
valid in a neighborhood of $M$. We mention here that we could
have worked with a simpler decomposition corresponding to a ``product'' cusp 
metric, if we did not want to evaluate the constants $C_2,C_3$. 
 
\begin{p}\label{prop3} 
The de Rham differential
\[d:\cun(X,\Lambda^k X)\to \cun(X,\Lambda^{k+1}(X))\]
restricts to a cusp differential operator 
\[d:\cun(\oX,\Lambda^k(\ctX))\to \cun(\oX,\Lambda^{k+1}(\ctX)).\] 
Its normal operator in the decomposition \eqref{decom} is 
\[\cN(d)(\xi)= 
\begin{bmatrix} 
d^M-i\xi\alpha_0\wedge&d^M\alpha_0\wedge\\ 
i\xi&-(d^M-i\xi\alpha_0\wedge) 
\end{bmatrix}.
\] 
\end{p} 
\begin{proof} 
Let $\omega\in\cun(\oX,\Lambda\ctX)$ and decompose it 
using \eqref{decom} into 
$\omega=\omega_1+V^0\wedge\omega_2$, where $V_0\lrcorner\omega_{1,2}=0$. 
Trivialize the bundle $\Lambda^*(\oX)$ with the local basis 
$V^0,\ldots,V^{n-1}$ and write for $j=1,2$ using \eqref{vi} 
\begin{align*} 
d\omega_j=&\sum_{i=1}^{n-1}V^i\wedge V_i(\omega_j) +V^0\wedge V_0(\omega_j)\\
=&d^M\omega_j -\alpha\wedge V_0(\omega_j) +V^0\wedge V_0(\omega_j)\end{align*} 
Now $dV^0=d\alpha\equiv d^M\alpha$ modulo $\cI^2$ (recall that the ideal
$\cI$ maps to $0$ under the normal operator map) so the announced formula
holds.
\end{proof} 
 
From the definition, the normal operator map is linear and multiplicative.
Moreover, it is a star-morphism, i.e., it commutes with taking the (formal)
adjoint.
This last property needs some explanation since the volume form used 
for taking the adjoint on the boundary must be specified. 
 
\begin{lemma} \label{lema6}
Let $P\in\Diffc(\oX,E,F)$ be a cusp operator and $P^*$ its adjoint with
respect to $g_0$. Then $\cN(P^*)(\xi)$ is the adjoint of $\cN(P)(\xi)$ 
with respect to the metric on $E_{|M},F_{|M}$ induced by restriction,
and the volume form ${a_0}^{1/2} dh_0$.
\end{lemma} 
\begin{proof} 
Since $\cN$ commutes with products and sums, 
it is enough to prove the result for a set of local generators of $\Diffc$.
Work in local coordinates $y_i$ on $M$. Then the vector fields 
$x^2\px$, $\partial_{y_i}$ and the set of smooth 
bundle morphisms $\{A\in\cun(\oX,E,F)\}$
form local generators for $\Diffc(\oX,E,F)$. 
From the definition of $\cN$ it follows that 
\begin{align*} 
\cN(x^2\px)(\xi)=&i\xi;\\ 
\cN(\partial_{y_i})(\xi)=&{\partial_{y_i}}_{|M};\\ 
\cN(A)(\xi)=&A_{|M}. 
\end{align*} 
Clearly, restriction of $A$ to $M$ commutes with taking the adjoint. 
The volume form of $g_0$ is
$dg_0=\sqrt{a\det(h_{ij})}x^2dx dy$. Since the derivative 
of $a$ and of the metric coefficients $h_{ij}$ with respect to 
$x^{2}\px$ belong to $\cI^2$, we see that $x^2\px^*\equiv-x^2\px+\cI^2$,
so
\[\cN((x^2\px)^*)(\xi)=\cN(-x^2\px)(\xi)=-i\xi=(i\xi)^*=(\cN(x^2\px)(\xi))^*.\]
Similarly,
\[(\partial_{y_i})^*=-\left(a\det(h_{ij})\right)^{-1/2}\partial_{y_i} 
(a\det(h_{ij}))^{1/2}\] 
restricts to $M$ to the adjoint of $\partial_{y_i}$ on $M$ 
with respect to the volume form $\sqrt{a\det(h_{ij})}_{|M}dy
={a_0}^{1/2} dh_0$. 
\end{proof} 
 
A related result is the conjugation invariance by
powers of $x$. Namely, if $P\in\Diffc$ and $s\in\cz$ then
$x^s Px^{-s}\in\Diffc$, and 
\[\cN(x^s Px^{-s})=\cN(P).\] 
 
The principal symbol of a cusp operator on $X$ extends as a map on the cusp
cotangent bundle down to $x=0$. In particular a cusp operator of positive
order cannot be elliptic. A cusp operator is called (cusp) elliptic if its 
principal symbol is invertible on $\ctsX\setminus\{0\}$ down to $x=0$.
In the rest of the article we refer to this property as ellipticity, as no
confusion can arise. A cusp operator is 
called \emph{fully elliptic} if it is elliptic and moreover its normal
operator is invertible for all values of $\xi\in\rz$.
 
\begin{te}\label{th4} 
The Laplacian $\Delta_{k,p}$ belongs to $x^{-2p}\Diffc^2(\oX,\Lambda^k(\ctX))$.
Moreover, $x^{2p}\Delta_{k,p}$ is fully elliptic if and only if the de Rham
coho\-mo\-lo\-gy groups $H^k_{\mathrm{dR}}(M)$ and $H^{k-1}_{\mathrm{dR}}(M)$
both vanish.
\end{te} 
\begin{proof} 
The principal symbol of the Laplacian 
of $g_p$ on $\Lambda^k X$ is simply $g_p$ times the identity. Since 
$x^{-2p}g_p=g_0$ (on the cotangent bundle!), it follows that 
$x^{2p}\Delta_{k,p}$ is elliptic down to $x=0$. 
 
Let $\delta_0^k, \delta_p^k$ be the formal adjoint of $d:\Lambda^k X 
\to\Lambda^{k+1}(X)$ with respect to $g_0$,
resp.\ $g_p$. Then $\delta_p^k=x^{(2k-n)p+2}\delta_0 x^{(n-2(k+1))p-2}$. 
What matters is that the sum of the exponents of $x$ equals $-2p$. 
By conjugation invariance, 
\[\cN(x^{2p}(d\delta_p+\delta_p d))=\cN(d\delta_0+\delta_0d).\] 
 
By Hodge theory (essentially, the fact that the range of an elliptic
operator is closed), the kernel of $\cN(d\delta_0+\delta_0d)(\xi)$ 
is isomorphic to the cohomology of the complex $(\Lambda^*(\ctX)_{|M}, 
\cN(d)(\xi))$.
 
We write $\cN(d)(\xi)=A(\xi)+B(\xi)$ where 
\begin{align*} 
A(\xi)=&\begin{bmatrix}0&0\\i\xi&0\end{bmatrix},& 
B(\xi)=&\begin{bmatrix} 
d^M-i\xi\alpha_0\wedge&d^M\alpha_0\wedge\\ 
0&-(d^M-i\xi\alpha_0\wedge) 
\end{bmatrix}.
\end{align*} 
We claim that for $\xi\neq 0$ the cohomology of $\cN(d)(\xi)$ vanishes.
The idea is to use again Hodge theory but with respect to the volume
form $dh_0$ on $M$. Then $B(\xi)^*$ anticommutes with $A(\xi)$ and similarly
$A(\xi)^*B(\xi)+B(\xi)A(\xi)^*=0$. Therefore
\[\cN(d)(\xi)\cN(d)(\xi)^*+ 
\cN(d)(\xi)^*\cN(d)(\xi)=\xi^2 I+B(\xi)^*B(\xi)+B(\xi)B(\xi)^*\] 
where $I$ is the $2\times 2$ identity matrix. So for $\xi\neq 0$
the Laplacian of $\cN(d)(\xi)$ is a strictly positive elliptic operator, hence 
it is invertible. 
 
Let us turn to the case $\xi=0$. We claim that the cohomology 
of $(\Lambda^*M\oplus\Lambda^{*-1}M, \cN(d)(0))$ is isomorphic to 
$H^*_{\mathrm{dR}}(M)\oplus H^{*-1}_{\mathrm{dR}}(M)$. Indeed, it
is enough to notice that 
\[\cN(d)(0)=\begin{bmatrix}1&-\alpha_0\wedge\\0&1\end{bmatrix} 
\begin{bmatrix}d^M&0\\0&-d^M\end{bmatrix} 
\begin{bmatrix}1&\alpha_0\wedge\\0&1\end{bmatrix}.\] 
In other words, the differential $\cN(d)(0)$ is conjugated to 
the diagonal de Rham differential, so they have isomorphic cohomology. 
\end{proof} 
 
\section{Proof of Theorem \ref{th3}} 
\begin{proof}[Proof of Theorem \ref{th3}] 
By Theorem \ref{th4}, the hypothesis implies that $\Delta_{k,p}$
is fully elliptic. By \cite[Theorem 17]{wlom}, the spectrum of $\Delta_{k,p}$
is pure-point and accumulates towards infinity according to \eqref{rg}, modulo
identification of the correct coefficients. To identify these coefficients
we use \cite[Proposition 14]{wlom} and Delange's theorem 
(\cite[Lemma 16]{wlom}. 
Note that $\sigma_2(\Delta_{k,p})$ is identically $1$ on the cosphere bundle.
The dimension of the form bundle equals $\dbinom{n}{k}$. This gives $C_1$ 
and $C_2$ using the above cited results. If we assume that $g_0$ is exact 
(i.e., the $1$-form $\alpha_0$ is exact and the function $a_0$ equals $1$)
then using a different function $x$ inside the cusp structure 
we can assume that $\alpha_0=0$. Thus for exact metrics the computation of
$\cN(d)$ gives 
\begin{equation}\label{cndkp}
\cN(x^{2p}\Delta_{k,p})(\xi)=\begin{bmatrix}\xi^2+\Delta_k^M&0\\0& 
\xi^2+\Delta_{k-1}^M 
\end{bmatrix}.
\end{equation}
This allows us to compute the integral from \cite[Proposition 14]{wlom} 
in terms of the zeta functions of the form Laplacians on $M$ with 
respect to $h_0$. 
\end{proof} 
 
\section{The magnetic Laplacian}\label{sectionmagn} 
\subsection{The magnetic Laplacian on a Riemannian manifold} 
A \emph{magnetic field} $B$ on the Riemannian manifold $(X,g)$
is simply an exact real valued $2$-form. A \emph{magnetic potential} 
$A$ associated to $B$ is a $1$-form such that $dA=B$. We
form the \emph{magnetic Laplacian} acting on $\cun(X)$: 
\[\Delta_A:=\dA^*\dA.\] 
This formula makes sense for complex-valued $1$-forms $A$.
Note that when $A$ is real, $\dA$ is a metric connection on the trivial 
bundle $\underline{\cz}$ with the canonical metric, and $\Delta_A$ 
is the connection Laplacian.

If we alter $A$ by adding to it a real exact form,
say $A'=A+dw$, the resulting magnetic Laplacian satisfies 
\[\Delta_{A'}=e^{-iw}\Delta_Ae^{iw}\] 
so it is unitarily equivalent to $\Delta_A$ in $L^2(X,g)$.
Therefore if $H^1_{\mathrm{dR}}(X)=0$ (for instance if $\pi_1(X)$ is finite)
then $\Delta_A$ depends, up to unitary equivalence, only on the
magnetic field $B$. This property is called \emph{gauge invariance}. 

In the literature one usually encounters gauge invariance as a
consequence of $1$-connectedness (i.e., $\pi_1=0$). But in dimensions 
at least $4$, every finitely generated group (in particular, every
finite group) can be realized as $\pi_1$ of a compact manifold. Thus
the hypothesis $\pi_1=0$ is unnecessarily restrictive.
 
While the properties of $\Delta_A$ in  $\rz^n$ with the flat metric
are quite well understood, the (absence of) essential spectrum of 
magnetic Laplacians on other manifolds has not been much studied so far.
One exception is the case of bounded geometry, studied in \cite{KS}. 
However our manifolds are not of bounded geometry because the injectivity
radius tends to $0$ at infinity.
 
\subsection{The cusp magnetic Laplacian } 
 Let $X$ be like in Section \ref{lcm}. 
We assume that the metric is of the form \eqref{pme} near $\partial X$, 
and moreover that $A\in\cun(\oX,\ctsX)$, in other words $A$ extends to 
a smooth section of $\ctsX$. We have seen in the introduction 
that the scalar Laplacian
(corresponding to $A=0$) always has continuous spectrum, at least
for the warped product metrics \eqref{mc}. We will show in Section
\ref{nonem} that this fact continues to hold for general
conformally cusp metrics.

Nevertheless, for most other $A$, $\Delta_A$ has pure-point spectrum and we can
describe its eigenvalue asymptotics. Note that $A$ below can be complex-valued.
\begin{te}\label{thmag} 
Let $p>0$, $g_p$ a metric on $X$ given by \eqref{pme} near $\partial X$ 
and $A\in\cun(X,T^*X)$ a complex-valued $1$-form satisfying near $\partial X$ 
\begin{equation}\label{aft} 
A=\varphi(x)\frac{dx}{x^2}+\theta(x)\end{equation} 
where $\varphi\in\cun(\oX)$ and $\theta\in\cun([0,\varepsilon)\times M, 
\Lambda^1(M))$. 
Assume that, on each connected component $M_0$ 
of $M$, either $\varphi_0:=\varphi(0)$ is not constant, or that
$\theta_0:=\theta(0)$ is not closed, or the cohomology class 
$[\theta_0]\in H^1_{\mathrm{dR}}(M_0)$ does not belong to the image of
\[2\pi H^1(M_0,\zz)\to H^1(M_0,\cz)\simeq H^1_{\mathrm{dR}}(M_0)\otimes\cz.\] 
Then $\Delta_A$ has pure-point spectrum. Its counting function satisfies
\begin{equation*} 
N_{A,p}(\lambda) \approx \begin{cases} 
C_1\lambda^{n/2}& 
\text{for $1/n< p<\infty$,}\\
C_2\lambda^{n/2}\log \lambda &\text{for $p=1/n$,}\\ 
C_3\lambda^{1/2p}&\text{for $0<p<1/n$} 
\end{cases} 
\end{equation*} 
in the limit $\lambda\to\infty$, where $C_1,C_2$ are given by
\eqref{c12} with $k=0$. 
\end{te} 
In particular, we note that $C_1,C_2$ do not depend on $A$. The hypotheses on
$A$ is another way of saying that $A$ is a smooth cusp $1$-form on $\oX$. 
\begin{proof} 
As before, we show that $x^{2p}\Delta_A$ is fully elliptic. To follow
the proof, the reader may as well assume that $\alpha_0=0$; the price to 
pay for this simplification is losing the information on $C_2$. 
We write $A=\varphi\cdot(dx/x^2+\alpha)+(\theta-\varphi\alpha)$ in order 
to use the decomposition \eqref{decom}. Using Proposition \ref{prop3} on 
functions, we get 
\[\cN(\dA)(\xi)= 
\begin{bmatrix} d^M-i\xi\alpha_0+i(\theta_0-\varphi_0\alpha_0)\\ 
i(\xi+\varphi_0)\end{bmatrix}.\] 
Assume that there exists $u\in\ker(\cN(x^{2p}\Delta_A)(\xi))$ different from
$0$. By elliptic regularity $u$ is smooth. 
We replace $M$ by one of its connected components on which $u$ does not
vanish identically, so we can suppose that $M$ is connected. 
Using Lemma \ref{lema6}, by integration by parts 
with respect to the volume form ${a_0}^{1/2}dh_0$ on $M$ 
and the metric $h_0$ on $\Lambda^1(M)$, we see
that $u\in\ker(\cN(\Delta_A)(\xi))$ implies $u\in\ker(\cN(\dA)(\xi))$. 
Set $\omega:=\theta_0-\alpha_0(\varphi_0+\xi)\in\Lambda^1(M)$. Then 
\begin{equation}\label{xvp}\begin{split}
(\xi+\varphi_0)u= &0,\\
(d^M+i\omega)u= &0,
\end{split}\end{equation} 
so $u$ is a global parallel section in the trivial bundle $\ucz$
over $M$, with respect to the connection $d^M+i\omega$. 
This implies 
\begin{equation*}\label{zudg} 
0=(d^M)^2 u=d^M(-iu\omega)=-i(d^M u)\wedge\omega-iud^M\omega=-iud^M\omega.
\end{equation*}  By uniqueness of 
solutions to ordinary differential equations, $u$ is never $0$, so 
$d^M\omega=0$. 
Furthermore, from Eq.\ \eqref{xvp}, we see that $\varphi_0$ equals the constant
function $-\xi$. Therefore $\omega=\theta_0$. In conclusion $\theta_0$ is closed 
and $\varphi_0$ is constant. It remains to show the assertion about the
cohomology class $[\theta_0]$. 
 
Let $\tM$ be the universal cover of $M$. Denote by ${\tilde{u}}$, ${\tilde{\theta}_0}$ the 
lifts of $u,\theta_0$ to $\tM$. The identity $(d^M+i\theta_0)u=0$ lifts to 
\begin{equation}\label{ecu} 
(d^{\tM}+i{\tilde{\theta}_0}){\tilde{u}}=0. 
\end{equation} 
The $1$-form ${\tilde{\theta}_0}$ is closed on the simply connected manifold 
$\tM$, hence it is exact 
(recall that by the universal coefficients formula, 
$H^1(\tM,\cz)=H_1(\tM,\cz)=H_1(\tM,\zz)\otimes\cz$, and $H_1(\tM,\zz)$ vanishes
as it is the abelianization of $\pi_1(\tM)$). Let $v\in\cun(\tM)$ 
be a primitive of ${\tilde{\theta}_0}$, i.e., $d^{\tM}v={\tilde{\theta}_0}$. Then ${\tilde{u}}$ can
be explicitly computed from \eqref{ecu}:
\[{\tilde{u}}=Ce^{-iv}\] 
for some constant $C\neq 0$.
 
The fundamental group $\pi_1(M)$ acts to the right on $\tM$ via
deck transformations. The condition that ${\tilde{u}}$ be the lift of 
$u$ from $M$ is the invariance under the action of $\pi_1(M)$, in other words 
\[{\tilde{u}}(y)={\tilde{u}}(y[\gamma])\] 
for all closed loops $\gamma$ in $M$. This is obviously equivalent to 
\begin{align*}
v(y[\gamma])-v(y)\in 2\pi\zz,&&\forall y\in \tilde{M}.\end{align*} 
Let $\tilde{\gamma}$ be the lift of $\gamma$ starting in $y$. Then
\begin{align*} 
v(y[\gamma])-v(y)=&\int_{\tilde{\gamma}}d^{\tM}v 
=\int_{\tilde{\gamma}}{\tilde{\theta}_0}=\int_\gamma \theta_0. 
\end{align*} 
Thus the solution ${\tilde{u}}$ is $\pi_1(M)$-invariant if and only if the cocycle $\theta_0$ 
evaluates to an integer multiple of $2\pi$ on each closed loop $\gamma$.
These loops span $H_1(M,\zz)$, so $[\theta_0]$ lives in the image of
$H^1(M,\zz)$ inside $H^1(M,\cz)=\Hom(H_1(M,\zz),\cz)$. Therefore the solution $u$ 
cannot be different from $0$ unless $\varphi_0$ is constant, $\theta_0$ 
is closed and $[\theta_0]\in 2\pi H^1(M,\zz)$.
 
We proved that under the hypotheses of the theorem, $\Delta_A$, which
belongs to $x^{-2p}\Diffc^2(M)$, is fully elliptic. Then by
\cite[Theorem 17]{wlom} we get the advertised eigenvalue asymptotics.
The constants $C_1, C_2$ are identified exactly as in the proof of
Theorem \ref{th3} (note that from \cite[Proposition 14]{wlom}, they depend
only on the principal symbol, which is unaffected by $A$).
\end{proof} 

Conversely, if $\varphi_0$ is constant, $\theta_0$ is closed and
$[\theta_0]\in 2\pi H^1(M,\zz)$ then ${\tilde{u}}=e^{-iv}$ like above is $\pi_1(M)$-invariant, so 
it is the lift of some $u\in\cun(M)$ which belongs to 
$\ker(\cN(x^{2p}\Delta_A)(\xi))$ for $-\xi$ equal to the constant value of $\varphi$. 
Thus $\Delta_A$ is not fully elliptic, hence it is not Fredholm between cusp
Sobolev spaces. Unfortunately, this fact alone does not allow us to 
conclude anything about its self-adjoint extension(s) in $L^2$, and more 
work is needed. In case of exact cusp metrics and real magnetic potentials
we give a converse to Theorem \ref{thmag} (and also to Theorem \ref{th:0})
in Section \ref{nonem}.
 
We have also had difficulties in extending the above proof to the magnetic 
Laplacian on $k$-forms.

Note that the hypotheses of Theorem \ref{thmag} are independent of the choice of
the magnetic potential $A$ in the class of cusp $1$-forms. 
Indeed, assume that $A'=A+dw$ for some $w\in\cun(X)$ is again a cusp $1$-form.
Then $dw$ must be itself a cusp form, so 
\[dw=x^2\px w \frac{dx}{x^2}+d^M w\in\cun(\oX,\ctsX).\] 
Write this as $dw=\varphi'\frac{dx}{x^2}+\theta'_x$. 
For each $x>0$ the form $\theta'_x$ is exact. By the Hodge decomposition theorem, 
the space of exact forms on $M$ is closed, 
so the limit $\theta'_0$ is also exact. Now $dw$ is an exact cusp form,
in particular it is closed. This implies that $x^2\px\theta'=d^M \varphi'$. 
Setting $x=0$ we deduce $d^M(\varphi')_{x=0}=0$, or equivalently
$\varphi'_{|x=0}$ is constant. Hence the conditions from Theorem \ref{thmag}
on the vector potential are satisfied simultaneously by $A$ and $A'$. 
 
\subsection{On the non-stability of the essential spectrum} 
In $\rz^n$, with the flat metric, it is shown in \cite{MP} that the essential 
spectrum of a magnetic Laplacian (for a bounded magnetic field)
is determined by the restriction of the magnetic field to the boundary of
a suitable compactification of $\rz^n$. In other words, 
only the behavior of the magnetic field at infinity plays
a r\^ole in the computation of the essential spectrum.
 
Moreover, \cite[Theorem 4.1]{KS} states that the non-emptiness
of the essential spectrum is preserved by the addition of a bounded magnetic 
field, again for the Euclidean metric on $\rz^n$. 
 
In contrast with these examples, Theorem \ref{thmag} indicates that in general 
the essential spectrum is not stable (and may even vanish) 
under compact perturbation of the magnetic field. 
This phenomenon is quite remarkable.

\subsection{The case $H^1(X)\neq 0$} The gauge invariance property 
does not hold in this case. Then the operators $\Delta_A$ and 
$\Delta_{A'}$ may have different essential spectra even when 
$d(A-A')=0$. We give below an example where the essential spectrum of
$\Delta_{A'}$ is empty, while that of $\Delta_{A}$ is not.

\begin{example}
Consider the manifold $\oX=(S^1)^{n-1}\times [0,1]$ with a metric $g_p$ as in
\eqref{pme} near the two boundary components. Let $\theta_j$ be variables
on the torus $(S^1)^{n-1}$, so $e^{i\theta_j}\in S^1$. Choose the vector
potentials $A$ to be $0$, and $A'$ to be the closed form $\mu d\theta_1$
for some $\mu\in\rz$. In this example the magnetic field $B$ vanishes.
We have mentioned in the introduction that the Laplacian on functions
has non-empty essential spectrum. However, the class $[i_M^*(A')]$ is
an integer multiple of $2\pi$ if and only if $\mu\in\zz$. Therefore by 
Theorem \ref{thmag}, for $\mu\in\rz\setminus\zz$ the essential spectrum 
of $\Delta_{A'}$ is empty. 
\end{example}
 
This effect is known as the Aharonov-Bohm effect \cite{AB}. The 
phenomenon we point out is surprising in light of what is known from 
the case of $\rz^n$ with an bounded obstacle. Indeed, the most you can expect 
in this case is a phase difference of the waves due to the choice of
two non-homotopic paths that circumvent the obstacle. Note that the
essential spectrum will not depend on the choice of $A$ in the case. 

\subsection{The case $H^1(X)= 0$}
Assume that $H^1(X,\rz)=0$ so $\Delta_A$ is
determined by $B$ up to unitary equivalence (by gauge invariance). 
Assume that $B$ is a smooth $2$-form on $\oX$.
which is closed on $X$ and such that its restriction to $X$ is exact. 
Then there exists $A\in \cun(\oX,\Lambda^1(\oX))$ such that $B=dA$
(since the cohomology of the de Rham complex on $\oX$ equals the
singular cohomology of $\oX$, hence that of $X$).
Assume moreover that the pull-back of $B$ to $M$ vanishes.
Then $B$ defines a relative cohomology class in $H^2(\oX,M,\rz)$. Besides,
the vanishing of $i_M^*(B)$ means that $i^*_M(A)$ is closed. It is easy
to see that $B$ determines a $1$-cohomology class $[i_M^*(A)]$ on $M$. 
Notice that $\partial[i^*_M(A)]=[dA]=[B]$, where
$\partial:H^1(M,\rz)\to H^2(\oX,M,\rz)$ is the connecting homomorphism
in the cohomology
long exact sequence of the pair $(\oX,M)$ with real coefficients. Consider
the same map but for integer coefficients, and the natural arrows from
$\zz$- to $\rz$-cohomology, which are just tensoring with $\rz$.
\begin{equation}\label{cohoz}
\begin{CD}
H^{1}(M,\rz)@>{\partial}>>H^2(X,M,\rz)\\
 @AA{\otimes\rz}A@AA{\otimes\rz}A\\
H^{1}(M,\zz)@>{\partial}>> H^2(X,M,\zz)
\end{CD}
\end{equation}

We see that $[B]$ lives in the image of $\otimes\rz$ if
$[A]$ does (conversely, if we assume additionally that $H_1(X,\zz)=0$,
which is the case for instance if $\pi_1(X)=0$, then $[B]$ is integer
if and only if $[A]$ is; this holds because $H^2(X,\zz)$ is torsion-free).
Thus from Theorem \ref{thmag} we get the following
\begin{cor}
Assume that the first Betti number of $X$ vanishes. Let $B$ be a smooth
exact $2$-form on $\oX$.
Assume that either $i_M^* B$ does not vanish identically, or that $i_M^* B=0$
and the cohomology class $[B]\in H^2(\oX,M,\rz)$ is not an integer multiple of
$2\pi$ (in the sense described above). Then the essential spectrum of
the associated magnetic Laplacian with respect to the metric \eqref{pme} is
empty.
\end{cor}

\begin{example}
Suppose that $B=df\wedge dx/x^2$ where $f$ is
a function on $X$ smooth down to the boundary $M$ of $X$. Assume that $f$
is not constant on $M$. Then the essential spectrum of the magnetic
operator (which is well-defined by $B$ if $H^1(X)=0$) is empty.
\end{example}
Note that in this example the pull-back to the border of
the magnetic field is null.
For such magnetic fields in flat $\rz^n$, 
the essential spectrum would be $[0,\infty)$.

\begin{example}\label{contreex} 
There exist Riemannian manifolds $X$ with gauge invariance
(i.e., $H^1(X,\rz)=0$) and a smooth
magnetic field $B$ with compact support in $X$ such
that the essential spectrum of the magnetic Laplacian is empty.
We choose $X$ to be the interior of a smooth manifold with boundary.
Assume that $H^1(X)=0$, for instance $X$ is simply connected. 
Assume moreover that $H^1(M,\rz)\neq 0$ where $M$ is the boundary of $X$. 
Note that these two assumptions are impossible to fulfill simultaneously
for orientable manifolds in dimension $3$ (see Remark \ref{rem3}).

We construct $A$ like in (\ref{aft}) satisfying the hypotheses
of Theorem \ref{thmag}. We take $\varphi_0$ to be constant. Let
$\psi\in\cun([0, \varepsilon))$ be a cut-off function
such that $\psi(x)=0$ for $x\in[3\varepsilon/4,\varepsilon)$ and 
$\psi(x)=1$ for $x\in[0, \varepsilon/2)$. Since $H^1(M)\neq 0$,
there exists a closed $1$-form $\beta$ on $M$ which is not exact.
Up to multiplying $\beta$ by a real constant, we can assume that the
cohomology class $[\beta]\in H^1_{\mathrm{dR}}(M_0)$ does not belong
to the image of
$2\pi H^1(M_0,\zz)\to H^1(M_0,\rz)\simeq H^1_{\mathrm{dR}}(M_0)$. Choose
$\theta_x$ to be $\psi(x)\beta$ for $\varepsilon>x>0$
and extend it by $0$ to $X$.

By Theorem \ref{thmag}, $\Delta_A$ has pure point spectrum. 
On the other hand, the magnetic field $B=dA=(d\psi/dx)\wedge\beta$
has compact support in $X$.
\end{example}

\begin{example}\label{ex12}
The simplest particular case of Example \ref{contreex}  is $X=\rz^2$,
but of course not with its Euclidean metric.
The metric \eqref{pme} could be for instance (in polar coordinates)
\[r^{-2p}(dr^2+d\sigma^2).\]
Here $M$ is the circle at infinity and $x=1/r$ for large $r$.

The product of this manifold with a closed, connected,
simply connected manifold $Y$ of dimension $k$ yields an example
of dimension $2+k$. Indeed, by the K\"unneth formula,
the first cohomology group of $\rz^2\times Y$ vanishes, while
$H^1(S^1\times Y)\simeq H^1(S^1)=\zz$.
\end{example}

Clearly $k$ cannot be
$1$ since the only closed manifold in dimension $1$ is the circle. Thus
dimension $3$ is actually exceptional.

\begin{ob}\label{rem3}
For orientable $X$ of dimension $3$ (the
most interesting for Physical purposes)
the assumptions $H^1(X)=0$ and $H^1(M)\neq 0$
cannot be simultaneously fulfilled.
Indeed, we have the following long exact sequence (valid actually
regardless of the dimension of $X$)
\[H^1(X)\to H^1(M)\to H^2(X,M).\]
If $\dim(X)=3$, the spaces $H^1(X)$ and $H^2(X,M)$ are dual by
Poincar\'e duality, hence $H^1(X)=0$ implies $H^2(X,M)=0$ and so
(by exactness) $H^1(M)=0$.
\end{ob}

\begin{example}
In the setting of Example \ref{ex12}, the relative cohomology group
$H^2(X,M)$ is one-dimensional. The integrality condition becomes
particularly nice in this case. Namely, let $B$ be a compactly
supported magnetic field (or more generally, a smooth exact $2$-form on 
$\oX$ vanishing at $M$) which is not zero in $H^2(X,M)$. Then 
$\lambda B$ is integral precisely for $\lambda$ in an infinite
discrete subgroup of $\rz$. Thus the multiples of $B$ for which the
magnetic Laplacian has non-empty essential spectrum are discrete and
periodic. The fact that there is essential spectrum for integral $[B]$
follows from Proposition \ref{thema}.  
\end{example}

\section{Schr\"odinger operators}\label{schro}

We now consider Schr\"odinger operators with magnetic 
potentials on $k$-forms.
An easy and well-known fact, which follows from Proposition
\ref{p:cut}, is the compactness of the resolvent of 
\eqref{eq:op} for any self-adjoint potential $V$ such that 
$V\rightarrow+\infty$ at infinity. Using the cusp algebra
we get much stronger results. We will apply our methods for 
$V$ belonging to $x^{-2p}\cun(\oX)$ but not necessarily positive. 
Unlike in the Euclidean case, for $p\geq 1/n$ the eigenvalue asymptotics do not
depend on $V$. The result holds moreover for endomorphism-valued 
potentials $V$.
 
\begin{p}\label{thschr}
Let $p>0$, $g_p$ the metric on $X$ given by \eqref{pme} near
$\partial X$ and $V\in x^{-2p}\cun(\oX, \End(\Lambda^k \oX))$ a self-adjoint 
potential. Set $V_0:=(x^{2p}V)_{|M}\in\cun(M, \End(\Lambda^k \oX)_{|M})$.
Assume that $V_0\geq 0$ (i.e., $V_0$ is semi-positive definite), and 
in each connected component of $M$ there exists
$z$ with $V_0(z)>0$. Let $A$ be any vector potential of the form \eqref{aft} near
$M=\partial X$.
Then the magnetic Schr\"odinger Laplacian 
\begin{equation}\label{eq:op} 
\Delta_{A,V}:=\dA^*\dA+\dA \dA^*+V
\end{equation} 
is essentially self-adjoint and has pure-point spectrum. 
Its eigenvalue counting function satisfies 
\begin{equation*} 
N_{A,V}(\lambda) \approx \begin{cases} 
C_1\lambda^{n/2}& 
\text{for $1/n< p<\infty$,}\\
C_2\lambda^{n/2}\log \lambda &\text{for $p=1/n$} \\
C_3 \lambda^{\frac{1}{2p}} &\text{for $p<1/n$} 
\end{cases} 
\end{equation*} 
with $C_1,C_2$ given by \eqref{c12}. 
\end{p} 
\begin{ob} 
Note that for $p\geq 1/n$, the asymptotic behavior of the eigenvalues 
depends neither on $V$ nor on $A$, i.e., $C_1,C_2$ depend only on the metric! 
\end{ob} 
 
\begin{proof} 
We claim that the operator $\Delta_{A,V}$ is fully elliptic.
Indeed, we write first $\cN(x^{2p}\Delta_{A,V})=\cN(x^{2p}\Delta_A)+V_0$.
For all $\xi\in\rz$, the operators $\cN(x^{2p}\Delta_A)(\xi)$ and $V_0$
are non-negative, so a solution of $\cN(x^{2p}\Delta_{A,V})(\xi)$
must satisfy
\begin{align*}
\cN(x^{2p}\Delta_A)(\xi)\phi=0,&& V_0\phi=0.
\end{align*}
By unique continuation, solutions of the elliptic operator 
$\cN(x^{2p}\Delta_A)(\xi)$ which are not identically zero 
on a given connected component of $M$ do not vanish anywhere on that component. 
But then $V_0\phi=0$ and the assumption on $V_0$ show that there are no 
non-zero solutions, in other words $\Delta_{A,V}$ is fully elliptic.

By \cite[Lemma 10 and Corollary 13]{wlom}, it follows that $\Delta_{A,V}$
is essentially self-adjoint with domain $x^{2p}H^2_c(M)$,
and has pure-point spectrum. Note that $\Delta_{A,V}$
is not necessarily positive since $V$ may take negative values. Write
$V=V^++V^-$, where $V^+_0=V_0$ while $V^-\in x^{1-2p}\cun(\oX)$.
Then $V^-$ is compact from $H^{2,2p}=\Dc(\Delta_{A,V^+})$ to $L^2$,
so $\Delta_{A,V}$ is a compact perturbation of a strictly positive
operator and therefore it is bounded below. 

By \cite[Proposition 14]{wlom},
for $p\geq 1/n$, the asymptotic behaviour of the function 
$N(\Delta_{A,V})$ as $\lambda\to\infty$ is 
given by \eqref{c12}, in particular it only depends on the principal symbol.

For $0<p<1/n$, again by \cite[Proposition 14]{wlom}, 
the coefficient $C_3$ in the spectral asymptotics 
of $\Delta_{A,V}$ only depends on $\cN(\Delta_{A,V})$.
\end{proof} 
 
The coefficient $C_3$ can be computed if we assume that the metric is exact 
and $A$ for instance vanishes. It then depends on the zeta function of 
$\Delta^M+V_0$. The computation is very similar to Theorem \ref{th3} 
and we leave it as an exercise to the interested reader. 

Specialize now to the case of functions. 
Then we can have pure-point spectrum and 
asymptotic laws for potentials which are unbounded below in some, but not all
directions towards infinity.
\begin{example}
Let $X$ be a Riemannian manifold with metric given by \eqref{pme} near infinity, 
and $p>0$.
Let $V\in x^{-2p}\cun(\oX)$ be such that $V_0=(x^{2p}V)_{|M}$ is non-negative 
and not identically zero on $M$. Then Proposition \ref{thschr} implies that
the Schr\"odinger operator $\Delta+V$ has pure-point spectrum and obeys
a generalized Weyl law. There exist such $V$ which tend to $-\infty$ towards
some part of $M$. For instance, assume $\supp(V_0)$ does not cover $M$; this
implies that near $M\setminus \supp(V_0)$, $V/x^{1-2p}$ has a limit, which 
can be very well negative. Now if $p>1/2$, then $x^{1-2p}$ goes to $\infty$
so $V$ tends to $-\infty$ as we approach $M\setminus \supp(V_0)$. 
\end{example}

\begin{ob}
This example is quite surprising in light of what happens in Euclidean
$\R^n$. Let $\overline{\R^n}$ be its spherical compactification, and 
$V\in \cC({\R^n};{\R})$ a potential which extends continuously to 
$\overline{\rz^n}$ with values in $\overline\rz$.  
Then $\sigma_{\rm ess}(\Delta +V)=\emptyset$ if and only if 
$V(x)=+\infty$ for $x\notin \R^n$. One can nicely see this in the setting of
$C^*$-algebras using the formalism from \cite{GI}. 
The same remains true for the Laplacian on a tree with potentials that extend
to the hyperbolic compactification of the tree (see \cite{G}), and 
for the natural generalization of these operators on a total Fock space 
\cite{GG2}. 
\end{ob}

\begin{example} Assume that the magnetic potential $A$
satisfies the hypotheses of Theorem \ref{thmag}. Then we can control
the essential spectrum of the magnetic Schr\"odinger operator
$\Delta_{A,V}$ even for potentials which tend to $-\infty$ at infinity
in \emph{all} directions. Namely, there exists a strictly positive constant $c$
depending on $g_p$ and $i^*_M A$, such that for any real potential 
$V\in x^{-2p} \cun(\oX)$ satisfying
\[|V_0|<c,\]
the operator $\Delta_{A,V}$ on functions has pure-point
spectrum and satisfies the conclusion of Proposition \ref{thschr}. 
Indeed, the constant
$c$ is the infimum over $\xi\in\rz$ of the first eigenvalue of
$\cN(x^{2p}\Delta_{A,V})(\xi)$, which can be seen to be strictly positive.
The proof is again the same as in Proposition \ref{thschr}: 
under the above hypotheses, the operator $\Delta_{A,V}$ is a fully 
elliptic cusp operator of order $(2,2p)$, so the conclusion follows 
by the above cited results of \cite{wlom}.
\end{example}

\begin{ob}
Consider a potential $V\in x^{-2p}\cun(\oX,\End(\Lambda^k(\oX))$ \emph{not
ne\-ces\-sa\-ri\-ly self-adjoint} (e.g., a possibly imaginary function for $k=0$).
Let $A$ be as in Proposition \ref{thschr}. Then the Schr\"odinger operator
$\Delta_{A,V}$ has pure-point spectrum, provided that the eigenvalues
of the limiting endomorphism $(x^{2p}V)_{|x=0}$ belong to
$\cz\setminus (-\infty,0]$ for all points on $M$. For functions, if 
$A$ satisfies the
hypotheses of Theorem \ref{thmag} (so $\Delta_A$ is fully elliptic) then
it is enough to ask that the eigenvalues of $(x^{2p}V)_{|x=0}$ do not
live in $(-\infty,-c)$ for some $c>0$. The Weyl-type
formulae do not hold however when $V$ is not self-adjoint.
\end{ob}

\section{Finite multiplicity of $L^2$ eigenvalues}\label{fm}

We will need in section \ref{nonem} a technical lemma that guarantees 
that the deficiency
indices of the Laplacian in the non-complete case are finite. In fact we 
can prove much more, namely that for an exact cusp metric and all $p>0$,
the multiplicity of every eigenvalue of 
the maximal extension of the Laplacian on forms and of the magnetic Laplacian
is finite. This may seem surprising since in the non-fully elliptic case the
operators are not Fredholm; on the other hand, this phenomenon was already 
noted in other Melrose-type algebras (see for instance \cite{melaps}).

\begin{lemma}\label{l:indices}
Assume that the metric $g_0$ is exact (see Definition \ref{defex}).
Then the multiplicity of all complex $L^2$ eigenvalues, in the sense of 
distributions, of the Laplacian $\Delta_{k,p}$ and of the magnetic Laplacian
$\Delta_A$ is finite. 
\end{lemma}
\begin{proof}
Let $\Delta$ denote one of our operators. We start by noticing that 
$\Delta$ can be seen as an unbounded operator in a larger $L^2$ space. 
Namely, let $L^2_{\varepsilon}$ be the completion of $\cunc(X,\Lambda^*X)$ 
with respect to the volume form $e^{-\frac{2\varepsilon}{x}}dg_p$ for some 
$\varepsilon>0$. Clearly then $L^2_{\varepsilon}$ contains $L^2$. The operator 
$\Delta$ is no longer symmetric in $L^2_{\varepsilon}$, however a distributional
solution of $\Delta-\lambda$ in $L^2$ is evidently also a distributional 
solution of $\Delta-\lambda$ in $L^2_{\varepsilon}$. Thus the conclusion 
will follow by showing that $\Delta-\lambda$ has in $L^2_{\varepsilon}$
a unique closed extension, which moreover is Fredholm. The strategy for this 
is by now clear. First we conjugate $\Delta-\lambda$ through the isometry 
\begin{align*}
L^2_{\varepsilon}\to L^2,&&\phi\mapsto e^{-\frac{\varepsilon}{x}}\phi.
\end{align*}
We get an unbounded operator $e^{-\frac{\varepsilon}{x}}\Delta e^{\frac{\varepsilon}{x}}$
in $L^2$ which is unitarily equivalent to $\Delta$ (acting in $L^2_{\varepsilon}$).
Essentially from the definition, 
\[\cN(x^{2p}e^{-\frac{\varepsilon}{x}}(\Delta-\lambda)e^{\frac{\varepsilon}{x}})(\xi)=
\cN(x^{2p}\Delta)(\xi+i\varepsilon)\]
(the term containing $\lambda$ vanishes under the normal operator since $p>0$). 
If the metric is exact, Eq.\ \eqref{cndkp} (for $\Delta=\Delta_{k,p}$) and a 
direct computation (for $\Delta=\Delta_A$) give 
\[\cN(x^{2p}\Delta)(\xi+i\varepsilon)=(\xi+i\varepsilon)^2+Q\]
for a certain explicit non-negative elliptic operator $Q$ on $M$. Since 
$\varepsilon\neq 0$, $\cN(x^{2p}\Delta)(\xi+i\varepsilon)$ 
is clearly invertible for $\xi\neq 0$. At $\xi=0$ 
we get $Q-\varepsilon^2$ which is also invertible, provided 
$\varepsilon^2\notin\Spec(Q)$. It is enough to choose $\varepsilon^2$ outside 
the spectrum of $Q$, which is pure-point since $M$ is closed. For such 
$\varepsilon$ the operator $\Delta-\lambda$ is unitarily equivalent to a 
fully elliptic cusp operator of order $(2,2p)$. It is then a general fact about 
the cusp algebra \cite[Theorem 17]{wlom}
that such an operator has a unique closed extension and has pure-point 
spectrum, in particular its eigenvalues have finite multiplicity.
\end{proof}

\section{Non emptiness of the essential spectrum in the
non fully-elliptic case}\label{nonem}

To complete our investigation it remains to compute the essential spectrum when
the operators are not fully-elliptic. We restrict ourselves to 
metrics which are of the form \eqref{mc} near $M$.

\begin{ob} \label{obwlog}
In the case of complete \emph{exact} cusp metrics, we can assume without loss
of generality that the metric is the ``toy metric'' \eqref{mc}. Indeed,
for an exact cusp metric (see Definition \ref{defex}),
by changing the function $x$ inside its cusp
structure, we can assume that $\alpha_{|M}=0$ (see \cite{melaps} and
\cite{wlom}). Moreover, using Theorem \ref{th:asymp} and Proposition
\ref{p:asymp}, if the metric is complete then, as far as the essential 
spectrum is concerned, we may replace
$h(x)$ in \eqref{cume2} by the metric $h_0:=h(0)$ on $M$, extended to a
symmetric $2$-tensor constant in $x$ near $M$.
\end{ob}

\subsection{The Laplacian acting on forms}
We will use the decomposition principle (Proposition \ref{p:cut}) to localize 
the computation of the essential spectrum to
the end $X':=(0,\varepsilon)\times M\subset X$. Set
\[\Krke:=L^2\left((0,\varepsilon), x^{(n-2k)p-2} dx\right).\]
Using \eqref{decom}, we get:  
\begin{equation}\label{decomph}  
L^2(X',\Lambda^k X)=\Krke\otimes 
\left(L^2(M,\Lambda^k M)\oplus L^2(M,\Lambda^{k-1} M)\right)
\end{equation}  
(the $k-1$ forms on $M$ appear as multiples of $dx/x^2$).
  
Now, we set $\Hr_0:=\Krke\otimes \ker(\Delta^M_k)$ and
$\Hr_1:=\Krke\otimes \ker(\Delta^M_{k-1})$. Then \eqref{decomph} becomes
\begin{equation}\label{dp} 
L^2(X',\Lambda^k X)=\Hr_0\oplus\Hr_1\oplus \Hr_h
\end{equation}
where (following \cite{lott}) the space of ``high energy forms'' $\Hr_h$ 
is by definition the orthogonal complement of $\Hr_0\oplus\Hr_1$.
  
\begin{p}\label{p:free}  
The Laplacian acting on $k$-forms $\Delta_{k,p}$ stabilizes the  
decomposition \eqref{dp} of $L^2(X',\Lambda^k X)$. 
Let $\Delta_{k,p}^0$, $\Delta_{k,p}^1$ and $\Delta_{k,p}^h$ be the
Friedrichs extensions of the restrictions of $\Delta_{k,p}$ to these spaces, 
respectively. The essential spectrum of $\Delta_{k,p}$ is the superposition
of the essential spectra of these three operators. Moreover, 
$\Delta_{k,p}^h$ has compact resolvent, and
\begin{align*}    
\Delta^0_{k,p}=\big(D^*D+c_0^2x^{2-2p}\big)\otimes 1, & 
& \Delta^1_{k,p}=\big(D^*D+c_1^2x^{2-2p}\big)\otimes 1,
\end{align*}  
where  
\begin{align}\label{cou}
c_0:=((2k+2-n)p-1)/2,&&c_1:=((2k-2-n)p+1)/2
\end{align}
and  $D:=x^{2-p}\px-c_0x^{1-p}$ acts in $\Krke$.
\end{p}  

Note that the essential spectra of $\Delta^0_{k,p}$, $\Delta^1_{k,p}$ are 
computed in Proposition \ref{p:essDelta}.

\begin{proof}
Since the de Rham operator and the Hodge star stabilize \eqref{dp}, so 
does the Laplacian. We will show that $\siges(\Delta_{k,p})$ only arises 
from $\Hr_0$ and $\Hr_1$. For this we use the cusp algebra
to perturb $\Delta_{k,p}$ into a fully elliptic cusp operator. 

Let $P$ denote the orthogonal 
projection in $L^2(M,\Lambda^k M\oplus\Lambda^{k-1} M)$ onto the     
space $\ker(\Delta^M_k)\oplus \ker(\Delta^M_{k-1})$ of harmonic forms.
Choose a real Schwartz cut-off function $\psi\in\cS(\rz)$ with $\psi(-C)=1$. 
Assume moreover
that the Fourier transform $\hat{\psi}$ has compact support (the reason for 
this assumption will appear later). Then $\psi(\xi)P$ defines a suspended 
operator of order $-\infty$ (see e.g., \cite[Section 2]{kso}). From 
\eqref{cndkp} we see that
\[\cN(x^{2p}\Delta_{k,p})(\xi)+\psi^2(\xi)P\]    
is strictly positive, hence invertible for all $\xi\in\rz$.     
By the surjectivity of the normal operator, there exists
$R\in\Psi_c^{-\infty}(X,\Lambda^k X)$ such that in the decomposition
\eqref{decom} over $M$,   
\[\cN(R)(\xi)=\psi(\xi)P.\] 
Fix a cut-off function $\phi$ with compact support in $X$ (this function 
appears in the proof of Proposition \ref{p:cut}), and another cut-off function
$\eta$ on $\oX$ which is $1$ near $M$ and such that $\eta\phi=0$.
By multiplying $R$ to the left and to the right by $\eta$ we can assume that
$R\phi=\phi R=0$, without changing $\cN(R)$. In fact we can give an explicit
formula for the Schwartz kernel of $R$:
\begin{equation}\label{skr}
\kappa_R(x,x',z,z')=\eta(x)\hat{\psi}\left(\frac{x-x'}{x^2}\right)\eta(x')
\kappa_P(z,z')\end{equation}
where $\kappa_P$ is the Schwartz kernel of $P$ on $M^2$. We assumed that 
$\hat{\psi}$ has compact support, thus $R$ acts on $\cunc(X',\Lambda^k X)$.
Notice that $R$ preserves the decomposition \eqref{dp}, and acts by $0$
on the space $\Hr_h$ of high energy forms. Let
\begin{equation}\label{rp}
R_p:=x^{-p}R\in \Psi_c^{-\infty,p}(X,\Lambda^k X).
\end{equation}
Then $R_p^*R_p\in \Psi_c^{-\infty,2p}(X,\Lambda^k X)$
is symmetric on $\cunc(X,\Lambda^k X)$ with respect to $dg_p$. Moreover, by the 
above discussion, $\Delta_{k,p}+R_p^* R_p$ is fully elliptic, so 
by \cite[Theorem 17]{wlom}, it is essentially self-adjoint in $L^2(M,dg_p)$
with domain $x^{2p}H^2(M)$, and has pure-point spectrum.

We apply now the decomposition principle to the non-negative operator
$\Delta_{k,p}+R_p^* R_p$. By Remark \ref{remfin}, 
the spectrum of $\Delta_{k,p}+R_p^* R_p$ can be computed on 
$X'=(0,\varepsilon)\times M$. Thus the Friedrichs extension
of $\Delta_{k,p}+R_p^* R_p$ on $X'$ also has pure-point spectrum; 
there, $\Delta_{k,p}+R_p^*R_p$ preserves the decomposition \eqref{dp}, 
and, by the
construction of $R_p$, it acts on the high energy forms as $\Delta_{k,p}$. 
Thus, it follows that over $X'$, the Friedrichs extension
of $\Delta_{k,p}$ with domain $\cunc(X')$ has pure point spectrum on $\Hr_h$.
Proving this fact has been the reason for introducing the operator $R_p$.

Again by the decomposition principle (Proposition \ref{p:cut}), 
$\siges(\Delta_{k,p})$ can be computed on the end $X'$. But there $\Delta_{k,p}$
also preserves \eqref{dp}, and we have just seen that its Friedrichs extension
has pure-point spectrum on the high energy forms.
It follows that $\siges(\Delta_{k,p})$ equals the superposition
of the essential spectra of the Friedrichs extensions
of $\Delta^0_{k,p}$ and $\Delta^1_{k,p}$.
A straightforward computation (see also \cite{francesca1}) shows that     
\begin{align}    
\Delta^0_{k,p}=&-x^{(2k-n)p}x^2\px x^{(n-2k-2)p}x^2\px\otimes 1\label{dez}    
\\    
\intertext{acting in $\Hr_0$. Similarly,}    
\Delta^1_{k,p}=&-x^2\px x^{(2k-2-n)p}x^2\px x^{(n-2k)p}\otimes 1\label{deu}    
\end{align}    
acting in $\Hr_1$. Expanding $D^*D$ from the definition, we get
the announced expressions for $\Delta^j_{k,p}$.   
\end{proof}  

We now compute the essential spectrum of $D^*D$:

\begin{lemma}\label{l:essD}
Let $D=x^{2-p}\px-c_0x^{1-p}$ act in 
$L^2\left((0,\varepsilon), x^{(n-2k)p-2} dx\right)$.
For $p\leq 1$, the essential spectrum of the operator $D^*D$ is $[0,\infty)$.
\end{lemma}
\proof
We conjugate $D^*D$ through the unitary transformation
\begin{align*}
L^2\left( x^{(n-2k)p-2} dx\right)\to L^2\left( x^{p-2} dx\right)&&
\phi\mapsto x^\alpha\phi
\end{align*} 
where $\alpha=(n-2k-1)p/2$. Let 
\[L(x)=\begin{cases}
-\ln(x)& \text{for $p=1$},\\
-\frac{x^{p-1}}{p-1}& \text{for $p\l 1$}.
\end{cases}\]
Under the change of variables $z:=L(x)$, $D^* D$ becomes the Euclidean Laplacian 
on $L^2\left((c,\infty), dz\right)$ for a certain $c$, plus 
a potential that goes to $0$ at infinity. 
\qed

We are now able to compute the essential spectrum of $\Delta_{k,p}$, and we
obtain Theorem \ref{th1} as a corollary of the next
proposition. 
We introduce the following set of \emph{thresholds}:
\begin{align*}
&\text{for $p<1$, }\cT=\begin{cases}
\emptyset, & \text{if $h^k(M)=h^{k-1}(M)=0$}, \\
\{0\}, & \text{otherwise}.\\
\end{cases}\\
&\text{for $p=1$, }\cT=\big\{c_i^2\in\{c_0^2,c_1^2\}; h^{k-i}(M)\neq 0\big\}.\\
&\text{for $p>1$, }\cT=\emptyset
\end{align*}
where $c_0,c_1$ are defined in \eqref{cou}.

\begin{p}\label{p:essDelta}
Assume either that the metric $g_0$ is exact and $p\leq 1$, or that
$g_0$ is the metric \eqref{mc} and $p>1$. Then the essential spectrum of 
$\Delta_{k,p}$ is given by 
\begin{eqnarray*}
\sigma_{\rm ess}(\Delta_{k,p})=[\inf(\cT), \infty)
\end{eqnarray*}
where $\cT$ is the set of thresholds.
In the incomplete case (i.e., $p>1$), this holds for all self-adjoint 
extensions of the operator $\Delta_{k,p}$.
\end{p}

Note that Theorem \ref{th:0} does not follow from Proposition \ref{p:essDelta}
because in the former we do not assume that the metric is exact.

\proof
Proposition \ref{p:free} identifies $\siges(\Delta_{k,p})$ with the
superposition of the essential spectra of $\Delta^i_{k,p}$ on $X'$.
For $p<1$, $\siges(\Delta^i_{k,p})=\siges(D^*D)$ since the
potential part tends to $0$ at infinity. For $p=1$ the potential 
in $\Delta^j_{k,p}$ is constant equal to $c_j^2$.
Therefore Lemma \ref{l:essD} gives the result for $p\leq 1$.

Let $p\g 1$. By Lemma \ref{l:indices} and by the Krein formula, all
self-adjoint extensions have the same essential spectrum so 
it is enough to consider the Friedrichs extension of $\Delta_{k,p}$.
We now use Proposition \ref{p:free}. 
The operator $D^*D$ is non-negative the spectrum of $\Delta_{k,p}^j$
is contained in $[\varepsilon^{2-2p}c_j^2,\infty)$, for $j=0,1$.
By Proposition \ref{p:cut}, the essential spectrum 
does not depend on the choice of $\varepsilon$. Now we remark that
$p>1$ implies $c_0,c_1\neq 0$. Indeed, in both cases the equality would
imply that $1/p\in\zz$, which is impossible. Thus
by letting $\varepsilon \to \infty$ we conclude that the essential
spectrum is empty. 
\qed

Note that for $p>1$, the ``toy metric'' $g_p$ given in \eqref{mc}
is essentially of \emph{metric horn} type \cite{lepe}.

\subsection{The magnetic Laplacian}

To further simplify the assumptions, let $A$ be given by \eqref{aft}
and define ${A_0}:={\varphi_0}dx/x^2+\theta_0$. Since
$\Delta_{{A_0}}$ is a perturbation of $\Delta_A$ by a first-order
operator small at infinity, by ellipticity it follows that we can
compute the essential spectrum of $\Delta_A$ using the metric
${g}_p$ given near $M$ by Eq \eqref{mc}, and the vector potential ${A_0}$
instead of $A$.

We have shown that the essential spectrum of the magnetic Laplacian
$\Delta_A$ is empty unless ${\varphi_0}$ is constant, $\theta_0$ is closed
and the cohomology class $[\theta_0]\in H^1(M)$ is an integer multiple of $2\pi$.
One can guess that in this last case,
the essential spectrum is not empty. Let us prove that this is so.

\begin{p}\label{thema}
Let $(X,g_p)$ be a Riemannian manifold with an exact cusp metric,
and $A\in\cun(\oX,\ctsX)$ a vector potential given by \eqref{aft} near $M$. 
Assume that ${\varphi_0}$ is constant, $\theta_0$ is closed
and the cohomology class $[\theta_0]\in H^1(M)$ is an integer multiple of
$2\pi$. Then 
\[\siges(\Delta_A)=\begin{cases}
[0,\infty)&\text{for $p<1$},\\
\left[\left(\frac{n-1}{2}\right)^2,\infty\right)&\text{for $p=1$},\\
\emptyset &\text{for $p>1$}.
\end{cases}\]
\end{p}
\begin{proof}
As shown above, we can assume that near $M$ we have
\begin{align*}
g_p=&x^{2p}\left(\frac{dx^2}{x^4}+h_0\right),&A=&Cdx/x^2+\theta_0.
\end{align*}
We decompose the space of $1$-forms as the direct sum \eqref{decom},
where now $V^0=x^{-2}dx$. Recall that $\delta_M$
is the adjoint of $d^M$ with respect to $h_0$.
Set $\dtt:=d^M+i\theta_0\wedge$. We compute
\begin{equation}\label{diam}\begin{split}
\dA=&
\begin{bmatrix}
d^M+i\theta_0\wedge\\
x^2\px+iC
\end{bmatrix}\\
\dA^*=& x^{-np}\begin{bmatrix}\delta_M-i\theta_0\lrcorner&
-(x^2\px+iC)
\end{bmatrix}x^{(n-2)p}\\
\Delta_A=& x^{-2p}(\dtt^*\dtt-(x^2\px+iC)^2-(n-2)px
(x^2\px+iC))
\end{split}\end{equation}
Since $d^M\theta_0=0$, the operator $\dtt$ on
$\cun(M,\Lambda^*(M))$ defines a
differential complex (i.e., $\dtt^2=0$); thus by Hodge theory, we 
decompose orthogonally
\[\cun(M)=\ker(\dtt)\oplus \im(\dtt^*)\]
and so
\begin{equation}\label{deco}
\cunc(M\times (0,\varepsilon))=\cunc(0,\varepsilon)\otimes\ker(\dtt)
\oplus \cunc(0,\varepsilon)\otimes\im(\dtt^*)
\end{equation}
Functions in the second term of this decomposition are called 
\emph{high energy functions}.
It is clear from \eqref{diam} that $\Delta_A$ preserves the decomposition
\eqref{deco}. Following the proof of Proposition \ref{p:free}, we will show 
that $\siges(\Delta_A)$ only arises from the first space in this decomposition.

Let $P$ be the orthogonal projection on $\ker\dtt$ in $L^2(M)$.
Choose a Schwartz cut-off function $\psi$ with $\psi(-C)=1$, whose
Fourier transform $\hat{\psi}$ has compact support
Then $\psi(\xi)P$ defines a suspended 
operator of order $-\infty$. We claim that
\[\cN(x^{2p}\Delta_A)(\xi)+\psi^2(\xi)P\]
is invertible for all $\xi$; indeed, this is a sum of non-negative operators 
for all $\xi$, the first of which is positive for $\xi\neq -C$ while at 
$\xi=-C$ we get $\dtt^*\dtt+P$ which is clearly invertible.

Define $R$ by Eq. \eqref{skr}
where now $\kappa_P$ is the Schwartz kernel of the above projector 
$P$ in $L^2(M)$, and consider the operator $R_p\in \Psi_c^{-\infty,p}(X)$ 
defined by \eqref{rp}. We assumed that 
$\hat{\psi}$ has compact support, thus $R$ acts on $\cunc(X)$.
Notice that $R$ preserves the decomposition \eqref{deco}, and acts by $0$
on the space of high energy functions. Then $\Delta_A+R_p^* R_p$ is fully elliptic
and symmetric, so by \cite[Theorem 17]{wlom} 
it is essentially self-adjoint with pure-point spectrum.

By the decomposition principle (Remark \ref{remfin}), the essential spectrum of 
$\Delta_A+R_p^* R_p$ can be computed on $X'=(0,\varepsilon)\times M$. 
Thus $\Delta_A+R_p^* R_p$ on $X'$ has also pure-point spectrum. Since $R_p$ is
$0$ on the high energy functions, it follows that the Friedrichs extension
of $\Delta_A$ with domain $\cunc(X')$ has pure point spectrum on 
the second factor of \eqref{deco}.

This is the place to mention that the hypotheses on $\theta_0$
imply that the kernel of $\dtt$
is $1$-dimensional. Indeed, let $v$ be a primitive of $\theta_0$ on the
universal cover of $M$. Then on different sheets of the cover, $v$ changes by
$2\pi\zz$ so $e^{iv}$ is a well-defined function on $M$. It is not hard to see
that this function spans $\ker(\dtt)$.

Again by the decomposition principle (Proposition \ref{p:cut}), 
$\siges(\Delta_A)$ can be computed on the end $X'$. But there $\Delta_A$
also preserves \eqref{deco}, and we have just seen that its Friedrichs extension
has pure-point spectrum on the high energy functions.
It follows that $\siges(\Delta_A)$ equals the essential spectrum of
the ordinary differential operator
\[-x^{-2p}((x^2\px+iC)^2+(n-2)px (x^2\px+iC))\]
in $L^2((0,\varepsilon),x^{np}dx)$. By the unitary transformation
$u\mapsto e^{iC/x}u$, we can suppose that $C=0$, i.e., we reduce to 
the case without magnetic field. Therefore
$\siges(\Delta_A)$ equals the essential spectrum 
of the scalar second-order differential operator
\[-x^{-2p}((x^2\px)^2+(n-2)px^3\px).\]
in $L^2((0,\varepsilon),x^{np-2}dx)$.
This is a particular case of the computation
for the Laplacian on forms. Namely, it is 
exactly the operator from \eqref{dez} for $k=0$, whose essential spectrum was
computed in Proposition \ref{p:essDelta}.
\end{proof}

\section{The nature of the essential spectrum} 

We consider the complete metric \eqref{mc} for $p\leq 1$.
The properties of the continuous spectrum of the Laplacian acting on 
functions are well known in this case. We now describe the nature of 
the essential spectrum of Laplacian acting on $k$-forms. 
 
For this metric, using Proposition \ref{p:free}, the study of the operator 
$\Delta_{k,p}$ is reduced to the analysis of \eqref{dez} and 
\eqref{deu}. 
We have already shown that the essential spectrum of $\Delta_{k,p}$ 
is given by $[\inf (\cT), \infty)$ where $\cT$ is the set of thresholds.

To pursue the analysis, one may use the so-called Mourre 
estimate introduced by E.~Mourre in \cite{mourre}. We refer to 
\cite{ABG} for a clear exposition of the theory. The major consequence 
of this estimate is to give a limit absorption principle. This allows 
one to deduce deep spectral results. 
 
By a decomposition principle, the study of the Laplacian can be reduced 
to the operators \eqref{dez} and \eqref{deu} on the half-line. 
The Mourre analysis of these operators is well understood. 
We refer to \cite{FH} for some general results. 
We deduce immediately the following 
corollaries: 
\begin{itemize} 
\item The singular continuous spectrum is empty. 
\item For each open set $J\subset\rz$ that does not contain a threshold, 
the number of eigenvalues is finite and of finite multiplicity.
\item The eigenvalues can accumulate only towards a threshold or towards 
infinity.
\item The multiplicity of the absolutely continuous spectrum is 
$h^k(M)1_{[c_0^2,\infty)}+h^{k-1}(M)1_{[c_1^2,\infty)}$.
\end{itemize} 

Note that we proved in Section \ref{fm} that all eigenvalues (including 
possibly the thresholds) have finite multiplicity. 

Compared to \cite{francesca2, francesca1}, our approach
based on the cusp calculus avoided
the analysis of the system of ordinary differential
equations which appears in the general case (the operator $\Delta^p_{M3}$
from loc.\ cit.). This system could in principle allow $\{0\}$ in 
the essential spectrum, and renders difficult
the analysis of absolutely/singular continuous spectrum. For instance,
for the metric (\ref{mc}), i.e. in the case $a\leq -1$ and
$b>0$ in the notation of \cite{francesca2,francesca1}, we can show that
the singular continuous spectrum $\sigma_{\rm{sc}}(\Delta_{k,p})$
is empty. We also do not have to restrict
ourself to the case $\partial X= S^{n-1}$. 
Finally, we compute the essential spectrum even for incomplete metrics
$g_p$, but we stress once more that we do not cover all the cases studied in 
loc.\ cit.

\appendix
\section{On the stability of the essential spectrum} 
 
We review some conditions under which the
essential spectrum of the Laplacian is stable by perturbation of 
the metric and by cutting a compact part of the manifold (the 
so-called decomposition principle).

\subsection{Perturbation of the metric}

We recall some general results from \cite[Section 5]{GG}. 
We stress that we suppose $X$ to be complete.
 
Let $X$ be a non-compact $\cC^1$ manifold, and $\mu$ the density of
a Riemannian metric $g$ on $X$. 
Let $\Hr:=L^2(X,\mu)$, and let $\Kr$ be the completion
of the space of continuous sections of $T^*X$ with compact support
under the natural norm
\[\|v\|_\Kr^2=\int_X\|v(x)\|_x^2 d\mu(x).\]

Let $d:\cC^1_c(X)\rightarrow\Kr$ be a closable first order operator.
It could be the de Rham differential, or
a magnetic de Rham differential, i.e. $d_{dR}+iA\wedge$ where $A$ is
a $\cC^1$ $1$-form, for instance. We keep the
notation $d$ for its closure. Its domain $\Gr$ is
the natural first order Sobolev space $\Hr^1_0$ defined in this
context as the closure of $\cC^1_c(X)$ under the norm
$$
\|u\|^2_{\Hr^1}=\int_X\Big(|u(x)|^2+\|d u(x)\|^2_x\Big)d\mu(x).
$$
Suppose that the injection from $\Hr^1_0(\Or)$ to $\Hr(\Or)$ is
compact, for all open bounded subsets $\Or\subset X$. The Laplacian 
of $d$ with respect to $g$ is given by $d^*d$.
 
\begin{te}\label{th:asymp}
Let $X$ be a $C^1$ non-compact manifold endowed with
a locally measurable metric $g_1$. Let $g_2$ be another locally
measurable metric, such that
\[\alpha(x)g_2(x)\leq g_1(x) \leq \beta(x)g_2(x),\]
with $\lim_{x\rightarrow\infty}\alpha(x)=\lim_{x\rightarrow\infty}
\beta(x)=1$. Suppose $X$ to be complete for one (thus also for the other) 
metric. Then the essential spectra of the Laplacians of $d$ with respect to
$g_1$ and $g_2$ (acting on functions) are the same.
\end{te} 
 
If one supposes moreover that $\mu(X)=\infty$, the convergence of 
$\alpha$ and $\beta$ can be assumed to hold in a weaker sense. 
 
We now give a version for the Laplacian acting on forms suitable for
our purpose in the context of cusp metrics. 

\begin{p}\label{p:asymp}
Let $X$ be a the interior of a compact smooth manifold. Let 
$g_1$ and $g_2$ be two cusp metrics on it.
Suppose that:
\begin{equation*}
\alpha(x)g_2(x)\leq g_1(x) \leq \beta(x)g_2(x)
\end{equation*}
such that $\lim_{x\rightarrow\infty}\alpha(x)=\lim_{x\rightarrow\infty}
\beta(x)=1$. Suppose $X$ to be complete for one (thus for the other) metric.
Then the essential spectra of the Laplacians acting on forms,
with respect to the two metrics are the same.
\end{p} 

The proof follows directly from the proof given at the end of
\cite{GG}[Section 5] in the case of forms. Indeed, in the
notation of \cite{GG}, for a complete cusp metric one checks easily that
$a^{\pm 1}\Gr\subset\Gr$. Note that the hypotheses
of this proposition can be significantly weakened.

\subsection{Removing a compact set} 
It is well-known that the essential spectrum of an elliptic differential 
operator on a complete manifold can be computed by cutting out a compact 
part and studying 
the Dirichlet extension of the remaining operator on the non-compact part 
(see e.g., \cite{donelili}). This result is obvious using Zhislin 
sequences, but the approach from loc.\ cit.\ fails in the non-complete case. 
For completeness, we give below a proof which has the advantage
to hold in a wider context (for cusp pseudodifferential operators).
 
We start with a general and easy lemma. 
We recall that a Weyl sequence for a couple $(H,\lambda)$ with $H$ a 
self-adjoint operator and $\lambda\in\rz$, is a sequence 
$\varphi_n\in\Dc(H)$ such that $\|\varphi_n\|=1$, 
$\varphi_n\rightharpoonup 0$ (weakly) and such that 
$(H-\lambda)\varphi_n\rightarrow 0$, as $n$ goes to infinity. 
The importance of this notion comes from the fact that 
$\lambda\in\sigma_{\rm ess}(H)$ if and only if there is Weyl sequence 
for $(H,\lambda)$. 
 
\begin{lemma}\label{l:cut} 
Let $H$ be a self-adjoint operator in a Hilbert space $\Hr$. 
Let $\varphi_n$ be a Weyl sequence for the couple $(H,\lambda)$. 
Suppose that there is a closed operator $\Phi$ in $\Hr$ such that: 
\begin{enumerate} 
\item $\Phi \Dc(H)\subset\Dc(H)$,\item $\Phi(H+i)^{-1}$ is compact, 
\item $[H,\Phi]$ extends to a bounded operator in $\Bc(\Dc(H),\Hr)$,
and $[H,\Phi](H+i)^{-1}$ is compact. 
\end{enumerate} 
Then there exists $\widetilde{\varphi}_n\in\Dc(H)$ such that 
$(1-\Phi)\widetilde{\varphi}_n$ is a Weyl sequence for $(H,\lambda)$.
\end{lemma} 
\proof 
First we note that $\Phi\varphi_n$ goes to $0$. Indeed, we have 
\begin{eqnarray*} 
\Phi\varphi_n&=&\Phi(H+i)^{-1}(H+i)\varphi_n\\ 
&=&\Phi(H+i)^{-1}\big((H-\lambda)\varphi_n +(i+\lambda) \varphi_n \big), 
\end{eqnarray*} 
the bracket goes weakly to $0$ and $\Phi(H+i)^{-1}$ 
is compact by (2). Similarly, using (3) we get that 
$[H,\Phi]\varphi_n\rightarrow 0$ too. 
 
There is $N$ such that for $n\geq N$, $\|(1-\Phi)\varphi_n\|\geq 1/2$ then we 
set $\widetilde\varphi_n=\varphi_n/\|(1-\Phi)\varphi_n\|$. 
We directly see that $(1-\Phi)\widetilde\varphi_n\rightharpoonup 0$. It 
remains to show that 
$(H-\lambda)(1-\Phi)\widetilde\varphi_n\rightarrow 0$ 
but this follows the fact that $[H,\Phi]\varphi_n\rightarrow 0$ and that 
$\varphi_n$ is a Weyl sequence for $(H,\lambda)$. 
\qed 
 
Let us now state the decomposition principle. Let $X$ be non-compact 
smooth Riemannian
manifold, and $K$ a compact sub-manifold with border of the same
dimension. We endow $X':=X\setminus K$ with the Riemannian
structure induced from $X$. 

\begin{p}\label{p:cut} 
Consider the Friedrichs extension of the magnetic
Laplacians $\Delta, \Delta'$ defined on compactly supported
forms on $X$
and $X'$, respectively. Then $\sigma_{\rm ess}(\Delta)=
\sigma_{\rm ess}(\Delta')$. 
\end{p} 
\proof Denote by $\eth$ the
closure of the operator $\dA+\dA^*$ acting on $\cunc(X)$,
and by $\eth'$ the closure of the same differential
operator acting on $\cunc(X')$.
We view $L^2(X',\Lambda^* X')$ as a closed subspace of $L^2(X,\Lambda^* X)$.
Set $\Qr(f,g)= \langle \eth f,\eth g\rangle$. We recall that domain of
$\Delta$ is given by  
\[\Dc(\Delta)=\{f\in \Dc(\eth ) \mid \text{$g\mapsto \Qr(f, g)$ is continuous 
for the $L^2$ norm}\}.\]

Let $\phi\in\cunc(X)$ be a non-negative cut-off function
such that $\phi|_K=1$. Let $\tif= 1-\phi$.
We claim that $f\in\Dc(\Delta)$ implies that $\tif
f\in\Dc(\Delta')$. Note first that the commutator $[\eth,\phi]$ is a smooth 
endomorphism with compact support (and so are all its derivatives). 
Thus $f\in\Dc(\eth)$ implies that $\tif f\in\Dc(\eth ')$, and $\eth'(\tif f)
=[\eth,\tif]f+\tif\eth f$. Therefore
\begin{eqnarray*}
\langle \eth u,\eth (\tif f)\rangle&=& \langle \eth u, [\eth ,\tif] f+ \tif \eth f\rangle =
\langle u, \eth ^* ([\eth ,\tif] f)\rangle + \langle \tif \eth u, \eth f\rangle\\
&=& \langle u, \eth ^* ([\eth ,\tif] f)\rangle + \langle \eth  (\tif u), \eth f\rangle
-\langle [\eth , \tif] u, \eth f\rangle
\end{eqnarray*}
which proves the claim. Similarly, we claim that $f\in\Dc(\Delta')$ implies 
$\tif f\in\Dc(\Delta)$. In fact, the previous argument shows that $\tif
f\in\Dc(\Delta')$. Now taking $\cchi$ another positive cut-off smooth
function such $\cchi=1$ on $K$ and $0$ on the support of $\tif$, for
$u\in \cunc(X)$ we notice that $\langle \eth u,\eth (\tif f)\rangle= \langle
\eth (1-\cchi)u,\eth (\tif f)\rangle$. This gives us the required continuity
and then that $\tif f\in\Dc(\Delta)$. Therefore if $f\in\Dc(\Delta)\cup
\Dc(\Delta')$, then $\tif f$ belongs to $\Dc(\Delta)\cap \Dc(\Delta')$. 
Moreover $\Delta \tif f=\Delta' \tif f$ since they are clearly 
both equal to $\Delta (\tif f)$ in the sense of distributions.

From the definition of the Friedrichs extension,
the domain of $\Delta$, $\Delta'$ is contained in $H^1_0(X)\cap H^2_\loc(X)$, 
respectively in $H^1_0(X')\cap H^2_\loc(X')$. Let $\Phi$ be the operator of 
multiplication by $\phi$. By the Rellich-Kondrakov lemma, the hypotheses of 
Lemma \ref{l:cut} are satisfied.  We apply Lemma \ref{l:cut} once for 
$H=\Delta$ and once for $H=\Delta'$ and since
$\widetilde\varphi_n\in\Dc(\Delta)\cup\Dc(\Delta')$, we have 
$\Delta (1-\phi)\widetilde \varphi_n= \Delta' (1-\phi)
\widetilde\varphi_n$. This proves the double inclusion of the essential spectra.
\qed

\begin{ob}\label{remfin}
More generally, consider a cusp operator of the form
$\tD=\Delta+R_p^*R_p$, where $R_p\in \Psi_c^{-\infty,p}(X)$ preserves
$\cunc(X')$ and is supported far from $K$ in the sense
\[\Phi R_p=R_p\Phi=0.\]
Such operators were explicitly constructed in the proofs of Propositions
\ref{p:free} and \ref{thema}. Note that 
$\tD$ preserves $\cunc(X')$, and that the domains of $\tD$, $\tD'$ are 
still contained in $H^1_0\cap H^2_\loc$. Then the proof of Proposition 
\ref{p:cut} holds for the Friedrichs 
extension of $\tD$. Namely if we set 
\[\Qr(f,g)= \langle \eth f,\eth g\rangle +\langle R_p f,R_p g\rangle,\]
the rest of the proof remains unchanged. 
\end{ob}

\bibliographystyle{plain} 

\begin{thebibliography}{1}
\bibitem{ABG} W.~Amrein, A.~Boutet de Monvel and V.~Georgescu, {\sl
  $C_{0}$-Groups, commutator methods and spectral theory of $N$-body
  Hamiltonians}, Birkh{\"a}user, Basel-Boston-Berlin, 1996. 

\bibitem{AB}
Y.~Aharonov and D.~Bohm,
{\sl Significance of electromagnetic potentials in the quantum theory, }
Phys. Rev. {\bf 115} (1959), no. 2, 485--491.

\bibitem{francesca2}
F.~Antoci,
{\sl On the absolute spectrum of the Laplace-Beltrami operator
for p-forms for a class of warped product metrics, }
preprint math.SP/0402391, to appear in Math. Nachrichten.
 
\bibitem{francesca1} 
F.~Antoci,
{\sl On the spectrum of the Laplace-Beltrami operator for $p$-forms
for a class of warped product metrics, } 
Adv. Math. {\bf 188} (2004), no. 2, 247--293. 
 
\bibitem{CFKS}
H.~Cycon, R.~Froese, W.~Kirsch and B.~Simon,
{\sl Schr\"odinger operators with application to quantum mechanics and global
geometry, } Texts and Monographs in Physics. Springer Study Edition.
Springer-Verlag, Berlin, 1987.

\bibitem{donelili} 
H.~Donnelly and P.~Li, 
{\sl Pure point spectrum and negative curvature for noncompact manifolds}, 
Duke Math. J. {\bf 46} (1979), no. 3, 497--503. 

\bibitem{FH}
R.~Froese and P.~Hislop, {\sl Spectral analysis of second-order
elliptic operators on noncompact manifolds, }
Duke Math. J. {\bf 58} (1989), no. 1, 103--129.

\bibitem{GG} 
V.~Georgescu and S.~Gol\'enia, {\sl Decay Preserving Operators
and stability of the essential spectrum, } preprint math.SP/0411489.

\bibitem{GG2} 
V.~Georgescu and S.~Gol\'enia, {\sl Isometries, Fock Spaces, and
  Spectral Analysis of Schr\"odinger Operators on Trees, } to appear in
J. Funct. Anal.

\bibitem{GI} 
V.~Georgescu and A.~Iftimovici,
{\sl Crossed products of $C\sp *$-algebras and spectral analysis of
  quantum  Hamiltonians, } Comm. Math. Phys. {\bf 228} (2002), no. 3, 519--560.

\bibitem{G} 
S.~Gol\'enia, {\sl  $C^*$-algebras of anisotropic Schr\"odinger
  operators on trees, } 
Ann. Henri Poincar\'e {\bf 5} (2004), no. 6, 1097--1115.

\bibitem{hm}
B.~Helffer and A.~Mohamed, {\sl Caract\'erisation du spectre essentiel
de l'op\'erateur de Schr\"odinger avec un champ magn\'etique, }
Ann. Inst. Fourier (Grenoble) {\bf 38} (1988), no. 2, 95--112.
 
\bibitem{KS} 
V.~Kondratiev and M.~Shubin {\sl Discreteness of spectrum
for the magnetic Schr\"odinger operators, } 
Comm. Partial Differential Equations {\bf 27} (2002), no. 3-4, 477--525. 

\bibitem{lani}
R.~Lauter and V.~Nistor,
{\sl On spectra of geometric operators on open manifolds and 
differentiable groupoids, }
Electron. Res. Announc. Amer. Math. Soc. {\bf 7} (2001), 45--53. 

\bibitem{lepe}
M.~Lesch and N.~Peyerimhoff, 
{\sl On index formulas for manifolds with metric horns, }
Comm. Partial Diff. Equ. {\bf 23} (1998), 649--684.

\bibitem{lott} J.~Lott, {\sl On the spectrum of a finite-volume
negatively-curved manifold, } Amer. J. Math. {\bf 123} (2001), no. 2, 
185--205.

\bibitem{MP} 
M.~M\u antoiu, R.~Purice and S.~Richard, {\sl Spectral and
Propagation Results for Magnetic Schr\"odinger Operators;
a $C^*$-Algebraic Framework,} preprint math.SP/0503046.

\bibitem{melaps}
R.~B.~Melrose,
{\sl The Atiyah-Patodi-Singer index theorem, }
Research Notes in Mathematics {\bf 4},
A. K. Peters, Wellesley, MA (1993).

\bibitem{meni96c}
R.~B.~Melrose and V.~Nistor,
{\sl Homology of pseudodifferential operators {I}. {M}anifolds with
boundary, }
preprint funct-an/9606005.

\bibitem{kso}
S.~Moroianu,
{\sl K-Theory of suspended pseudo-differential operators, }
$K$-Theory {\bf 28} (2003), 167-181.
 
\bibitem{wlom} 
S.~Moroianu,  
{\sl Weyl laws on open manifolds, } preprint 
math.DG/0310075.
 
\bibitem{mourre} 
E.~Mourre, {\sl Absence of singular continuous spectrum for 
certain self-adjoint operators, } Comm.\ Math.\ Phys.\ {\bf 91} 
(1981) 391--408.

\bibitem{shubin} 
M.~Shubin, {\sl Essential self-adjointness for semi-bounded magnetic  
Schr\"odinger operators on non-compact manifolds, }
J. Funct. Anal. {\bf 186} (2001), no. 1, 92--116. 
 
\bibitem{LT}  
M.~Loss and B.~Thaller, {\sl Scattering of particles by long-range 
  magnetic fields, }  Ann. Physics {\bf 176} (1987), no. 1, 159--180. 
 
\end{thebibliography}
 
\end{document}